\documentclass[12pt]{amsart}
\usepackage{amssymb}
\usepackage{amsbsy}
\usepackage{amscd}
\usepackage[mathscr]{eucal}
\usepackage{verbatim}
\usepackage{epic}
\usepackage{eepic}
\makeatletter
\def\cal{\mathcal}
\def\Bbb{\mathbb}
\def\frak{\mathfrak}

\newenvironment{pf*}[1]{\proof[#1]}{\endproof}
\renewcommand{\labelenumi}{(\theenumi)}%
\newcommand{\rom}{\textup}
\makeatother

\newenvironment{aenume}{%
  \begin{enumerate}%
  \renewcommand{\theenumi}{\alph{enumi}}%
  \renewcommand{\labelenumi}{\theenumi)}%
  }{\end{enumerate}}
\makeatletter
\@ifclasslater{amsart}{1999/11/24}{}{
\renewcommand*\subjclass[2][1991]{%
  \def\@subjclass{#2}%
  \@ifundefined{subjclassname@#1}{%
    \ClassWarning{\@classname}{Unknown edition (#1) of Mathematics
      Subject Classification; using '1991'.}%
  }{%
    \@xp\let\@xp\subjclassname\csname subjclassname@#1\endcsname
  }%
}
\renewcommand{\subjclassname}{%
  \textup{1991} Mathematics Subject Classification}
\@xp\let\csname subjclassname@1991\endcsname \subjclassname
\@namedef{subjclassname@2000}{%
  \textup{2000} Mathematics Subject Classification}
}
\makeatother
\newtheorem{Theorem}[equation]{Theorem}
\newtheorem{Corollary}[equation]{Corollary}
\newtheorem{Lemma}[equation]{Lemma}
\newtheorem{Proposition}[equation]{Proposition}

\theoremstyle{definition}

\newtheorem{Notation}[equation]{Notation}

\newtheorem{Conjecture}[equation]{Conjecture}

\theoremstyle{remark}
\newtheorem{Remark}[equation]{Remark}

\numberwithin{equation}{section}

\newcommand{\thmref}[1]{Theorem~\ref{#1}}
\newcommand{\secref}[1]{\S\ref{#1}}
\newcommand{\lemref}[1]{Lemma~\ref{#1}}
\newcommand{\propref}[1]{Proposition~\ref{#1}}

\newcommand{\lsp}[2]{\,{}^{#1}{#2}}

\newcommand{\C}{{\Bbb C}}
\newcommand{\Z}{{\Bbb Z}}
\newcommand{\Q}{{\Bbb Q}}
\newcommand{\R}{{\Bbb R}}
\newcommand{\proj}{{\Bbb P}}
\newcommand{\CP}{\proj}

\newcommand{\SU}{\operatorname{\rm SU}}
\newcommand{\GL}{\operatorname{GL}}

\newcommand{\U}{\operatorname{\rm U}}

\newcommand{\algsl}{\operatorname{\frak{sl}}} 

\newcommand{\End}{\operatorname{End}}
\newcommand{\Hom}{\operatorname{Hom}}
\newcommand{\Ext}{\operatorname{Ext}}
\newcommand{\Ker}{\operatorname{Ker}}

\newcommand{\Ima}{\operatorname{Im}}

\newcommand{\im}{\mathop{\text{\rm im}}\nolimits}
\newcommand{\rank}{\operatorname{rank}}

\newcommand{\tr}{\operatorname{tr}}

\newcommand{\pd}[2]{\frac{\partial#1}{\partial#2}}

\newcommand{\ve}{\varepsilon}

\newcommand{\linf}{\ell_\infty}
\newcommand{\shfO}{\mathcal O}
\newcommand{\dslash}{/\!\!/} 
\newcommand{\bp}{{\widehat\proj}^2}
\newcommand{\bM}{{\widehat M}}

\newcommand{\Supp}{\operatorname{Supp}}
\newcommand{\ch}{\operatorname{ch}}
\newcommand{\Wedge}{{\textstyle \bigwedge}}

\newcommand{\Fin}{F^{\text{\rm inst}}}
\newcommand{\Finz}{\mathcal F^{\text{\rm inst}}}
\newcommand{\q}{\mathfrak q}
\newcommand{\hT}{\widetilde T}

\begin{document}
\title[Instanton counting on blowup. I]
{Instanton counting on blowup. I.
\\ $4$-dimensional pure gauge theory}
\author{Hiraku Nakajima}
\address{Department of Mathematics, Kyoto University, Kyoto 606-8502,
Japan}
\email{nakajima@math.kyoto-u.ac.jp}
\thanks{The first author is supported by the Grant-in-aid
for Scientific Research (No.13640019, 15540023), JSPS}

\author{K\={o}ta Yoshioka}
\address{Department of Mathematics, Faculty of Science, Kobe University,
Kobe 657-8501, Japan}
\email{yoshioka@math.kobe-u.ac.jp}
\subjclass[2000]{Primary 14D21; Secondary 57R57, 81T13, 81T60}

\begin{abstract}
  We give a mathematically rigorous proof of Nekrasov's conjecture: the
  integration in the equivariant cohomology over the moduli spaces of
  instantons on $\mathbb R^4$ gives a deformation of the
  Seiberg-Witten prepotential for $N=2$ SUSY Yang-Mills theory.
  Through a study of moduli spaces on the blowup of $\mathbb R^4$, we
  derive a differential equation for the Nekrasov's partition
  function.  It is a deformation of the equation for the
  Seiberg-Witten prepotential, found by Losev et al., and further
  studied by Gorsky et al.
\end{abstract}

\maketitle

\section*{Introduction}




Let $M(r,n)$ be the framed moduli space of torsion free sheaves $E$ on
$\proj^2$ with rank $r$, $c_2 = n$, where the framing is a
trivialization of the restricition of $E$ at the line at infinity
$\linf$. There is a natural action of an $(r+2)$-dimensional torus
$\hT$, coming from the symmetry of the base space $\C^2 =
\proj^2\setminus\linf$ and the change of the framing.

Nekrasov's partition function \cite{Nek} is the generating function of
the integral of the equivariant cohomology class $1\in H^*_{\hT}(M(r,n))$:
\begin{equation*}
   Z(\ve_1,\ve_2,\vec{a};\q)
   = \sum_{n=0}^\infty \q^n \int_{M(r,n)} 1,
\end{equation*}
where $\ve_1$, $\ve_2$, $\vec{a} = (a_1,\dots,a_r)$ are generators of
$H^*_{\hT}(\mathrm{pt}) = S^*(\operatorname{Lie}\hT)$. When $n > 0$,
$M(r,n)$ is {\it noncompact\/} and the integration is given by {\it
  formally\/} applying the localization formula in the equivariant
cohomology. Then the integration $\int_{M(r,n)} 1$ is a rational
function in $\C(\ve_1,\ve_2,a_1,\dots,a_r)$. (A precise definition
will be given in the main body of the paper.)

Nekrasov conjectures that 
\(
      \Fin(\ve_1,\ve_2,\vec{a};\q)
       = \ve_1\ve_2\log Z(\ve_1,\ve_2,\vec{a};\q)
\)
is regular at $\ve_1, \ve_2 = 0$, and
\(
    \Fin(0,0,\vec{a};\q)
\)
is the instanton part of the Seiberg-Witten prepotential for $\mathcal
N=2$ supersymmetric $4$-dimensional gauge theory \cite{SW} with gauge
group $\SU(r)$.  Nekrasov's definition is mathematically rigorous. The
Seiberg-Witten prepotential is also rigorously defined by certain
period integrals of hyperelliptic curves (the so-called Seiberg-Witten
curves). The relation between the two, which is rather natural from a
physical point of view, can be considered as a mathematically well
formulated conjecture. It is very similar to the mirror symmetry. The
Nekrasov partition function is a counterpart of the Gromov-Witten
invariants and is on the `symplectic' side. Seiberg-Witten
prepotential is on the `complex' side.

Let us briefly recall the history on Donaldson invariants and
Seiberg-Witten prepotential. A reader can read the main body of the
paper without knowing the history, but then he/she loses the
motivation why we study Nekrasov's partition function.
In \cite{Wit} Witten described Donaldson invariants as the correlation
functions of certain operators in a twisted $\mathcal N=2$
supersymmetric Yang-Mills theory. Several years later Seiberg-Witten
found that the prepotential, which controls the physics of the theory,
can be computed via the periods of hyperelliptic curves \cite{SW}.
Then Moore-Witten studied Donaldson invariants using the
Seiberg-Witten prepotential \cite{MW}. In particular, they derived the
blowup formula for Donaldson invariants originally given by
Fintushel-Stern~\cite{Fintushel-Stern:1996}. These arguments were
physical and have no mathematically rigorous justification so far. It
was very misterious why Donaldson invariants are related to periods of
Seiberg-Witten curves.
Nekrasov's conjecture can be considered as a first step towards the
understanding of the misterious relation.

The main result in this paper can be summarized as follows. 
We consider a similar partition function defined via the framed moduli
space $\bM(r,k,n)$ on the blowup $\widehat{\C}^2$. We also introduce
an `operator' $\mu(C)$ associated with the exceptional set $C$.
We then show that the {\it correlation functions} 
\(
    \sum_{n=0}^\infty \q^n \int_{\bM(r,0,n)} \mu(C)^d
\)
vanish for $d=1,\dots, 2r-1$. This simplest case of the blowup formula
gives a differential equation \eqref{eq:ind} satisfied by
$Z(\ve_1,\ve_2,\vec{a};\q)$. We call it the {\it blowup equation}. The
blowup equation is a deformation of the differential equation
\eqref{eq:Wh2} for the Seiberg-Witten prepotential originally found in
the study of the contact term in the twisted $\mathcal N=2$
supersymmetric gauge theory by Losev et al.~\cite{LNS1,LNS2}. This
equation was derived also from the Seiberg-Witten curve in the
frame work of Whitham hierarchies by Gorsky et al.~\cite{GM3}. (A
self-contained proof will be given in \cite{lecture}.) By
Edelstein et al.~\cite{EMM} the equation determines the instanton
corrections recursively (see also \cite{Mar} and the references
therein).
An immediate application is an affirmative solution of Nekrasov's
conjecture: $\Fin(0,0,\vec{a};\q)$ is the instanton part of the
Seiberg-Witten prepotential.

Our strategy goes in the inverse direction of the above mentioned
history. We define the operator $\mu(C)$, mimicking the definition of
the similar operator for Donaldson invariants. Our vanishing is
well-known for Donaldson invariants (see e.g., \cite{FM}) and our
proof is exactly the same.  But this rather trivially looking
observation leads to the powerful blowup equation as we just
mentioned. (We eventually recover the whole Fintushel-Stern's formula
for arbitrary $d$ and its higher rank analog given in \cite{MaM} in
\secref{sec:blowupformula}.)  Let us remark that a relation between
Fintushel-Stern's blowup formula and the Whitham hierarchy was pointed
out in \cite[\S3]{LNS2}.
We also remark that there was an approach to Fintushel-Stern's blowup
formula based on Uhlenbeck (partial) compactifications of framed
moduli spaces \cite{Bryan:1997}.
The use of the simplest (or lowest) case of the blowup formula to
derive constraint is {\it not\/} a new idea in the context of
Donaldson invariants. The proof in \cite{Fintushel-Stern:1996} was
done essentially by this idea. G\"ottsche determined the wall-crossing
formula also by this idea \cite{Go}.

The paper is organized as follows. 
In \secref{sec:SWprep} we recall the Seiberg-Witten prepotential.
In \S\S\ref{sec:framed}, \ref{sec:mod_blowup}, we define framed moduli
spaces of coherent torsion free sheaves on the plane and its blowup.
We define an action of an $({r+2})$-dimensional torus $\hT$ on framed
moduli spaces, classify the fixed point set and determine the weights
of tangent spaces at fixed points. In \secref{sec:Hilbert} we consider
a natural $K$-theory analog of Nekrasov's partition function and
identify it with a Hilbert series of the coordinate ring of the framed
moduli spaces. This result partly explains why Nekrasov's partition
function is natural. But this reformulation is also used to prove the
simplest blowup formula.  \secref{sec:rank1} is a small detour. We
study the rank $1$ case, i.e., when the moduli spaces are Hilbert
schemes of points. Nekrasov's partition function and its blowup
formula is easy to derive, but some feature of the general cases can
be seen in this simplest case.  \S\S\ref{sec:counting},
\ref{sec:limit} are main part of this paper.  We introduce the
operator $\mu(C)$ and derive the blowup equation. We then prove
Nekrasov's conjecture. In \secref{sec:blowupformula} we derive the
full blowup formula for our correlation function of $\mu(C)$. In
\secref{sec:gauge} we consider the case when the gauge group is not
necessarily $\SU(r)$. Moduli spaces of torsion-free sheaves do not
have generalization to other gauge groups, so we are forced to use
Uhlenbeck (partial) compactifications. Our formulation in
\secref{sec:Hilbert} has a modification by using Uhlenbeck
compactifications. We then prove the blowup equation under some
technical assumptions on geometric properties of Uhlenbeck
compactifications.

In this paper, we treat only the pure gauge theory. Theories with
matters, as well as the inclusion of higher Casimir operators (i.e.,
we integrate more general cohomology classes other than $1$), will be
studied in the later series.

Our project started in 1997 together with I.~Grojnowski. The first
goal was a new proof of the blowup formula for Betti numbers of moduli
spaces originally given by the second author \cite{Yoshioka:1996}.
This part was finished soon afterward, and was reported by the first
author at Workshop on Complex Differential Geometry, 14-25 July 1997,
Warwick and at Verallgemeinerte Kac-Moody-Algebren, 19-25 July 1998,
Oberwolfach. (There are closely related results by W-P.~Li and Z.~Qin
\cite{LQ1,LQ2,LQ3}. We explain this result in \cite{lecture}.)
We then tried to give a new proof of Fintushel-Stern's blowup formula
for Donaldson invariants. The technique was to use the localization
theorem in the equivariant cohomology of the framed moduli space on
the blowup, which is basically the same technique taken in this paper.
But we did not understand how to take the `nonequivariant limit' since
a naive limit diverges. Thus we did not succeed at that time, and a
failure report was given by the first author at a workshop at RIMS
Kyoto, June 2000 \cite{blowup}.
The correct choice of limit is provided via the use of the Nekrasov's
partition function, and we finally succeed this time. And we get
Nekrasov's conjecture as a bonus.

While we were writing this paper, we were informed that Nekrasov and
Okounkov also proved Nekrasov's conjecture \cite{NO}. Their method is
totally different from ours.

After writing the first version of this paper, the authors gave series
of lectures on the subject at ``Workshop on algebraic structures and
moduli spaces'', July 14--20, 2003, Universite de Montreal. The reader
can find physical backgrounds and various related topics in the
lecture notes \cite{lecture}.

\subsection*{Acknowledgement}
The authors are grateful to I.~Grojnowski for discussion in the early
stage of our project.
They also thank the referees for helpful suggestions and comments.

\section{Seiberg-Witten prepotential}\label{sec:SWprep}

In this section, we briefly recall the definition of the
Seiberg-Witten prepotential for the sake of the reader. See
\cite[\S2]{lecture} for detail and proofs.

We consider a family of hyperelliptic curves parametrized by
$\vec{u} = (u_2,\dots,\linebreak[0]u_r)$:
\begin{equation*}
   C_{\vec{u}} : \Lambda^r\left(w + \frac1w\right)
   = P(z) = z^r + u_2 z^{r-2} + u_3 z^{r-3} + \cdots + u_r.
\end{equation*}
We call them {\it Seiberg-Witten curves}. The projection
$C_{\vec{u}}\ni (w,z)\mapsto z\in \proj^1$ gives a structure of
hyperelliptic curves. The parameter space $\{ \vec{u} \in \C^{r-1}\}$
is called the {\it $u$-plane}.

Let $z_1,\dots,z_r$ be the solutions of $P(z) = 0$. We will work on a
region of the $u$-plane where $|z_\alpha - z_\beta|$, $|z_\alpha|$ are
much larger than $|\Lambda|$.
We can find $z_\alpha^\pm$ near $z_\alpha$ such that $P(z_\alpha^\pm)
= \pm 2\Lambda^r$ as $|u| \gg |\Lambda|$. These are the $2r$-branched
points of the projection $C_{\vec{u}}\to\proj^1$.

The hyperelliptic curve $C_{\vec{u}}$ is made of two copies of the
Riemmann sphere, glued along the $r$-cuts between $z_\alpha^-$ and
$z_\alpha^+$ ($\alpha = 1,\dots, r$), as usual.
Let $A_\alpha$ be the cycle encircling the cut between $z_\alpha^-$
and $z_\alpha^+$. We have $\sum_\alpha A_\alpha = 0$. We take cycles
$B_\alpha$ so that $\{ A_\alpha, B_\alpha \mid \alpha=2,\dots, r\}$
form a symplectic basis of $H_1(C_{\vec{u}},\Z)$, i.e., $A_\alpha
\cdot A_\beta = 0 = B_\alpha\cdot B_\beta$, $A_\alpha\cdot B_\beta =
\delta_{\alpha\beta}$ for $\alpha,\beta=2,\dots, r$. (The cycle $A_1$
is omitted.) See \cite{lecture} for the precise choice.

Let us define the {\it Seiberg-Witten differential\/} by 
\begin{equation*}
   dS = - \frac1{2\pi 
                 }z \frac{dw}w
   = -\frac1{2\pi 
     }
     \frac{z P'(z) dz}{\sqrt{P(z)^2  - 4\Lambda^{2r}}}
   .
\end{equation*}
It is a meromorphic differential having poles at $\infty_\pm$. We
define functions $a_\alpha$, $a^D_\beta$ on the $u$-plane
by
\begin{equation*}\label{eq:def_a}
   a_\alpha = \int_{A_\alpha} dS, \qquad
   a^D_\beta = 2\pi \sqrt{-1} \int_{B_\beta} dS,
   \qquad \alpha=1,\dots,r,\ \beta=2,\dots,r.
\end{equation*}
We have $\sum_\alpha a_\alpha = 0$. In the gauge theory side, $\vec{a}
= (a_1,\dots,a_r)$ will be the coordinate system on $\operatorname{Lie}T$.

It can be shown that there exists a locally defined function $\mathcal
F(\vec{a};\Lambda)$ on the $\vec{u}$-plane such that
\begin{equation*}\label{eq:SWprep}
   a_\alpha^D = - \frac{\partial\mathcal F}{\partial a_\alpha}.
\end{equation*}
It is called the {\it Seiberg-Witten prepotential}. Note that
\begin{equation*}\label{eq:period}
   \tau_{\alpha\beta}
   = - \frac1{2\pi\sqrt{-1}}
     \frac{\partial^2 \mathcal F}{\partial a_\alpha \partial a_\beta}
\end{equation*}
is the period matrix of $C_{\vec{u}}$.

One can show that $\mathcal F$ has the following behaviour at $\Lambda\to 0$:
\begin{equation}\label{eq:pert}
  \mathcal F =
  \sum_{\alpha < \beta}\left[ (a_\alpha-a_\beta)^2
    \log\left(\frac{\sqrt{-1}(a_\alpha-a_\beta)}{\Lambda}\right)
    -\frac32 (a_\alpha-a_\beta)^2\right] 
+ \Lambda^{2r}\times O(\Lambda^{2r}).
\end{equation}
The first part (resp.\ second part $\Lambda^{2r}\times
O(\Lambda^{2r})$) is called the {\it perturbative part\/} (resp. {\it
  instanton part\/}) of the prepotential. For the choice of the branch
of $\log$, see \cite[\S2]{lecture}.

We use terminology for root systems of Lie algebras. The change is
useful for considering generalization to other gauge groups (see
\secref{sec:gauge}).

We consider $\vec{a}$ as an element of the Cartan subalgebra
$\mathfrak h$ of $\algsl_r$. Let $\Delta\subset\mathfrak h^*$ be the
set of roots. We take standard simple roots $\alpha_i\in \mathfrak
h^*$ and simple coroots $\alpha_i^\vee\in\mathfrak h$
($i=1,\dots,r-1$), i.e., $\alpha_i = (0,\dots, 0, \overset{i}{1},
\overset{i+1}{-1},0,\dots,0)$.
Let $\Delta_+$ denote the set of positive roots, i.e.,
\(
   \Delta_+ = \{ e_{\alpha,\beta} =
   (0,\dots, 0, \overset{\alpha}{1}, 0, \dots, 0,
   \overset{\beta}{-1},0,\dots,0) \mid \alpha < \beta \}.
\)
If $\vec{a} = (a_1,\dots, a_r)$, then $\langle \vec{a},
e_{\alpha,\beta} \rangle = a_\alpha - a_\beta$. We write $\vec{a} =
\sum_i a^i \alpha_i^\vee$.
Let $Q$ be the coroot lattice of $\mathfrak h$, i.e.,
\(
   Q = \{ \vec{k} = (k_1,\dots,k_r)\in \Z^r \mid \sum_\alpha k_\alpha = 0 \}.
\)
We write $\vec{k} = \sum k^i\alpha_i^\vee$ as above.

The perturbative part of the prepotential is rewritten as
\begin{equation*}
  \sum_{\alpha \in \Delta_+}
  \left[ \langle \vec{a},\alpha \rangle^2
    \log\left(\frac{\sqrt{-1}\langle \vec{a},\alpha \rangle}
   {\Lambda}\right)-\frac32 \langle \vec{a},\alpha \rangle^2\right].
\end{equation*}
The period matrix is
\begin{equation}\label{eq:tau}
     \tau_{ij} = -\frac{1}{2\pi \sqrt{-1}}
\frac{\partial^2 {\mathcal F}}{\partial a^i \partial a^j}
  = \frac{\sqrt{-1}}{\pi} \sum_{\alpha\in\Delta_+}
   \langle \alpha_i^\vee, \alpha\rangle 
   \langle \alpha_j^\vee, \alpha\rangle 
     \log\left(\frac{\sqrt{-1}\langle \vec{a},\alpha\rangle}
       {\Lambda}\right) + \Lambda^{2r}\times O(\Lambda^{2r}).
\end{equation}

In our proof of Nekrasov's conejcture, we use the following two equations:
\begin{align}
  \pd{{\mathcal F}}{\log \Lambda}&=-2r u_2,
  \label{eq:RG}
\\
  \pd{u_2}{\log \Lambda}=&-\frac{2r}{\pi \sqrt{-1}}
  \pd{u_2}{a^i}\pd{u_2}{a^j}\frac{\partial}{\partial\tau_{ij}}
  \log \Theta_E(\vec{0}|\tau),
  \label{eq:contactterm}
\end{align}
where
\begin{equation}\label{eq:Theta}
     \Theta_E(\vec{\xi} | \tau)
   = \sum_{\vec{k}\in Q} \exp\left(
     \pi\sqrt{-1} \sum_{i,j}\tau_{ij} k^i k^j
     + 2\pi\sqrt{-1}\sum_i k^i(\xi^i + \frac12)\right).
\end{equation}
The first equation \eqref{eq:RG} is called the {\it renormalization
  group equation}, and was obtained by \cite{STY}. (See also
\cite{Matone,EY,HKP1}.)

The second equation \eqref{eq:contactterm} is called the {\it contact
  term equation}. It was originally found in the context of the
$\mathcal N=2$ supersymmetric gauge theory \cite{LNS1,LNS2}, and
derived also from the above Seiberg-Witten curve in a mathematically
rigorous way \cite{GM3}.

\newcommand{\Phys}{\mathrm{Phys}}
\begin{Remark}
  In order to get the exact match with the physics literature, we need
  to note $\vec{a} = -\sqrt{-1}\vec{a}^{\Phys}$,
  $u_p = -u_p^{\Phys}$.
\end{Remark}

\section{Framed moduli spaces on the projective plane}\label{sec:framed}

In this section, we define framed moduli spaces on $\proj^2$ and study
their basic properties. All of results are straightforward
generalizations of the corresponding results for Hilbert schemes on
$\C^2$, which were explained in \cite{Lecture}. In fact, the results
were obtained long time ago and mentioned in
\cite[Exercise~5.15]{Lecture}.

Let $M(r,n)$ be the framed moduli space of torsion free sheaves on
$\proj^2$ with rank $r$ and $c_2 = n$, which parametrizes isomorphism
classes of $(E,\Phi)$ such that 
\begin{enumerate}
\item $E$ is a torsion free sheaf of $\rank E =r$, $\langle c_2(E),
[\proj^2]\rangle =n$ which is locally free in a neighborhood of
$\linf$,
\item $\Phi \colon E|_{\linf} {\overset{\sim}{\to}}
\shfO_{\linf}^{\oplus r}$ is an isomorphism called `framing at infinity'.
\end{enumerate}
Here $\linf=\{[0:z_1:z_2] \in \proj^2 \} \subset \proj^2$ is the line
at infinity.
Notice that the existence of a framing $\Phi$ implies $c_1(E)=0$.

The framed moduli spaces were constructed by
Huybrechts-Lehn~\cite{HL-frame} (in more general framework). The
tangent space is $\Ext^1(E,E(-\linf))$ and the obstruction space is
$\Ext^2(E,E(-\linf))$. In our situation, we have the following
vanishing theorem:
\begin{Proposition}\label{prop:vanish}
$\Hom(E,E(-\linf)) = \Ext^2(E,E(-\linf)) = 0$.
\end{Proposition}

\begin{proof}
By the Grothendieck-Serre duality theorem, $\Ext^2(E,E(-\linf))$ is
the dual of $\Hom(E,E(-2\linf))$. We shall show that
$\Hom(E,E(-k\linf)) = 0$ for any $k\in \Z_{> 0}$.

From a short exact sequence
\begin{equation*}
0 \to E(-(k+1)\linf) \xrightarrow{\text{mult. by $z_0$}} E(-k\linf) \to 
 E(-k\linf)\otimes {\shfO}_{\linf} \to 0,
\end{equation*}
we obtain an exact sequence
\begin{equation*}
0 \to \Hom(E, E(-(k+1)\linf)) \to \Hom(E, E(-k\linf)) 
  \to \Hom(E, E(-k\linf)\otimes\shfO_{\linf}).
\end{equation*}
Since the restriction of $E$ to $\linf$ is trivial, we have
\begin{equation*}
   \Hom(E, E(-k\linf)\otimes\shfO_{\linf}) = 0.
\end{equation*}
Hence we get
\begin{equation*}
   \Hom(E,E(-\linf)) \cong \Hom(E, E(-2\linf)) \cong \cdots
   \cong \Hom(E,E(-k\linf)) \cong \cdots.
\end{equation*}
But $\Hom(E,E(-k\linf))\cong \Ext^2(E,E((k-3)\linf))^*$ vanishes for
sufficient large $k$ by the Serre vanishing theorem. Thus we get the
assertion.
\end{proof}

\begin{Corollary}
$M(r,n)$ is a nonsingular variety of dimension $2nr$.
\end{Corollary}

\begin{proof}
  This follows from the above vanishing theorem together with the
  Riemann-Roch formula.
\end{proof}

In fact, we have another way to define the framed moduli space and
prove this corollary in our setting.
By a result of Barth~\cite{Barth} (see \cite[Theorem~2.1]{Lecture} for
the proof), we have an isomorphism between $M(r,n)$ and the quotient
space of $B_1,B_2 \in \End(\C^n)$, $i \in \Hom(\C^r,\C^n)$ and $j \in
\Hom(\C^n,\C^r)$ satisfying
\begin{enumerate}
\item $[B_1,B_2]+ij=0$,
\item there exists no proper subspace $S \subsetneq \C^n$ such that
$B_\alpha(S) \subset S$ \textup{(}$\alpha=1,2$\textup{)} and $\im i
\subset S$
\end{enumerate}
modulo the action of $\GL_n(\C)$ given by
\begin{equation*}
g \cdot (B_1,B_2,i,j)=( g B_1 g^{-1}, g B_2 g^{-1}, gi, jg^{-1}).
\end{equation*}
We say $(B_1,B_2,i,j)$ is {\it stable\/} when it satisfies the
condition~(2). It can be shown that the differential of the defining
equation~(1) is surjective and the action is free on stable
points. This shows the smoothness of $M(r,n)$. (See \cite[\S3]{Lecture}.)

Let $M_0(r,n)$ be the framed moduli space of ideal instantons on $S^4
= \C^2\cup\{\infty\}$, that is
\begin{equation*}
   M_0(r,n) = 
   \bigsqcup_{n'=0}^n M_0^{\operatorname{reg}}(r,n')\times S^{n-n'}\C^2,
\end{equation*}
where $M_0^{\operatorname{reg}}(r,n')$ is the framed moduli space of
genuine instantons on $S^4$ and $S^k\C^2$ is the $k$th symmetric
product of $\C^2$. We endow a topology to $M_0(r,n)$ as in
\cite[4.4]{DK}.
By a result of Donaldson~\cite{Don} (which is based on the ADHM
description \cite{ADHM}), $M_0^{\operatorname{reg}}(r,n)$ can be
identified with the framed moduli space of {\it locally free
sheaves\/} on $\proj^2$, and also with the open subset of the space of
linear data $(B_1,B_2,i,j)$ with an extra condition that the
transposes $\lsp{t}{B_1}$, $\lsp{t}{B_2}$, $\lsp{t}{j}$ satisfy the
above condition~(2). Then by \cite[3.4.10]{DK} together with
\cite[Chapter~3]{Lecture}, $M_0(r,n)$ can be identified (as a
topological space) with
\begin{equation}\label{eq:alg_quot}
 \left.\left\{ (B_1,B_2,i,j) \mid [B_1,B_2]+ij=0 \right\}\right.\dslash
 \GL_n(\C),
\end{equation}
where $\dslash$ denotes the affine algebro-geometric quotient. The
open locus $M_0^{\operatorname{reg}}(r,n)$ consists of closed orbits
$\GL_n(\C)\cdot(B_1,B_2,i,j)$ such that the stabilizer is trivial.

As in \cite[Chapter~3]{Lecture}, $M(r,n)$ has a structure of
hyper-K\"ahler manifold of dimension $4nr$ if we put the standard
inner products on $\C^n$ and $\C^r$.
In fact, $M(r,n)$ is isomorphic to the hyper-K\"ahler quotient
\begin{equation}
\label{eq:hyper}
\left\{ (B_1,B_2,i,j) \left|
\hspace{-0.3cm}
\begin{minipage}[m]{9cm}
\begin{enumerate}
\def\labelenumi{(\theenumi)}
\def\theenumi{\roman{enumi}}
\item $[B_1,B_2]+ij=0$
\item $[B_1, B_1^\dagger] + [B_2, B_2^\dagger] + ii^\dagger -
        j^\dagger j = \zeta\operatorname{id}$
\end{enumerate}
\end{minipage}
\right\}\right/ \U(n),
\end{equation}
where $(\ )^\dagger$ is the Hermitian adjoint and $\zeta$ is a fixed
positive real number. This hyper-K\"ahler structure plays no role
later.

By these descriptions via linear data, we have a projective morphism
\begin{equation*}
   \pi\colon M(r,n) \to M_0(r,n),
\end{equation*}
where we endow $M_0(r,n)$ with a scheme structure by the
description~\eqref{eq:alg_quot}. (See \cite[3.51]{Lecture}.) In terms
of the original definition as framed moduli spaces, the corresponding
map between closed points can be identified with
\begin{equation}\label{eq:map_pi}
   (E,\Phi) \longmapsto
   ((E^{\vee\vee},\Phi), \operatorname{Supp}(E^{\vee\vee}/E))\in
   M_0^{\operatorname{reg}}(r,n')\times S^{n-n'}\C^2.
\end{equation}
where $E^{\vee\vee}$ is the double dual of $E$ and
$\operatorname{Supp}(E^{\vee\vee}/E)$ is the support of
$E^{\vee\vee}/E$ counted with multiplicities. Note that $E^{\vee\vee}$
is a locally free sheaf. (This identification can be proved easily
from results in \cite[Chapters~2,3]{Lecture} and details were given in
\cite{VV}.)

\begin{Remark}
Morphisms from moduli spaces of semistable torsion-free sheaves to
moduli spaces of ideal instantons on general projective surfaces were
constructed by J.~Li \cite{Li} and Morgan \cite{Mor} in this way. (See
also \cite[\S8.2]{HL-book}.)
But it is not clear that the scheme structure is the same as one
given above.
\end{Remark}

Let $T$ be the maximal torus of $\GL_r(\C)$ consisting of diagonal
matrices and let $\hT = \C^*\times\C^* \times T$. We define an action
of $\hT$ on $M(r,n)$ as follows: For 
$(t_1,t_2)\in \C^*\times\C^*$, let $F_{t_1,t_2}$ be an automorphism of 
$\proj^2$ defined by
\[
    F_{t_1,t_2}([z_0: z_1 : z_2]) = [z_0: t_1 z_1 : t_2 z_2].
\]
For $\operatorname{diag}(e_1,\dots,e_r)\in T$ let $G_{e_1,\dots,e_r}$
denote the isomorphism of $\shfO_{\linf}^{\oplus r}$ given by
\[
    \shfO_{\linf}^{\oplus r}\ni (s_1,\dots, s_r) \longmapsto
     (e_1 s_1, \dots, e_r s_r).
\] 
Then for $(E,\Phi)\in M(r,n)$, we define
\begin{equation}\label{eq:action}
    (t_1,t_2,e_1,\dots,e_r)\cdot (E,\Phi)
    = \left((F_{t_1,t_2}^{-1})^* E, \Phi'\right),
\end{equation}
where $\Phi'$ is the composite of homomorphisms
\begin{equation*}
   (F_{t_1,t_2}^{-1})^* E|_{\linf} 
   \xrightarrow{(F_{t_1,t_2}^{-1})^*\Phi}
   (F_{t_1,t_2}^{-1})^* \shfO_{\linf}^{\oplus r}
   \longrightarrow \shfO_{\linf}^{\oplus r}
   \xrightarrow{G_{e_1,\dots, e_r}} \shfO_{\linf}^{\oplus r}.
\end{equation*}
Here the middle arrow is the homomorphism given by the action.

In a similar way, we have a $\hT$-action on $M_0(r,n)$. The map
$\pi\colon M(r,n)\to M_0(r,n)$ is equivariant.

\begin{Lemma}\label{lem:ADHMaction}
These actions can be identified with the actions on the linear data
defined by
\begin{equation*}
    (B_1, B_2, i, j) \longmapsto
     (t_1 B_1, t_2 B_2, i e^{-1}, t_1 t_2 e j), \qquad
   \text{for $t_1, t_2\in \C^*$,
             $e = \operatorname{diag}(e_1,\dots, e_r)\in(\C^*)^r$}.
\end{equation*}
Note that this action preserves the equation $[B_1,B_2]+ij=0$ and the
stability condition, and commutes with the action of
$\GL_n(\C)$. Hence it induces an action on $M(r,n)$ and $M_0(r,n)$.
\end{Lemma}

\begin{proof}
The sheaf $E$ is given as the middle cohomology group of the complex
\newcommand{\OCP}{{\shfO_{\proj^2}}}
$$
V \otimes \OCP(-1) \xrightarrow[a \, = \, \left(
 \begin{smallmatrix}z_0 B_1 -z_1 \\ z_0 B_2 -z_2 \\ z_0 j\end{smallmatrix}
\right)]{}
\begin{matrix} V \otimes \OCP \\
           \oplus      \\
         V \otimes \OCP  \\
           \oplus      \\
         W \otimes \OCP  \\
          \end{matrix}
\xrightarrow[
b \, = \, \left( \begin{smallmatrix}
-(z_0 B_2 - z_2) & z_0 B_1 - z_1 & z_0 i
\end{smallmatrix} \right)]{}
V \otimes  \OCP(1).
$$
(See \cite[Chapter~2]{Lecture}.)
Let us pull back this complex by $F_{t_1,t_2}^{-1}$:
$$
V \otimes \OCP(-1) \xrightarrow[a \, = \, \left(
 \begin{smallmatrix}z_0 B_1 -t_1^{-1}z_1 \\ z_0 B_2 -t_2^{-1}z_2 \\
   z_0 j
\end{smallmatrix}
\right)]{}
\begin{matrix} V \otimes \OCP \\
           \oplus      \\
         V \otimes \OCP  \\
           \oplus      \\
         W \otimes \OCP  \\
          \end{matrix}
\xrightarrow[
b \, = \, \left( \begin{smallmatrix}
-(z_0 B_2 - t_2^{-1} z_2) & z_0 B_1 - t_1^{-1} z_1 & z_0 i
\end{smallmatrix} \right)]{}
V \otimes  \OCP(1).
$$
Under the isomorphism
\begin{equation*}
\begin{matrix} V \otimes \OCP \\
           \oplus      \\
         V \otimes \OCP  \\
           \oplus      \\
         W \otimes \OCP  \\
          \end{matrix}
\ni
\begin{pmatrix}
  v_1 \\ v_2 \\ w
\end{pmatrix}
\longmapsto
\begin{pmatrix}
  t_2^{-1} v_1 \\ t_1^{-1} v_2 \\ w
\end{pmatrix},
\end{equation*}
the kernel of $b$ is mapped to the kernel of
\begin{equation*}
\left( \begin{matrix}
-(z_0 t_2 B_2 - z_2) & z_0 t_1 B_1 - z_1 & z_0 i
\end{matrix} \right).
\end{equation*}
Also under the above isomorphism, the image of $a$ is mapped to the image of
\begin{equation*}
\frac1{t_1t_2}
\begin{pmatrix}
z_0 t_1 B_1 -z_1 \\ z_0 t_2 B_2 - z_2 \\ z_0 t_1 t_2 j
\end{pmatrix}.
\end{equation*}
Thus the pull-back sheaf $(F_{t_1,t_2}^{-1})^* E$ corresponds to the data
$(t_1 B_1, t_2 B_2, i, t_1 t_2 j)$.
Composing the change of the framing by $G_{e_1,\dots, e_r}$, we get
the assertion.
\end{proof}


\begin{Proposition}\label{prop:fixedpoint}
\rom{(1)}
$(E,\Phi)\in M(r,n)$ is fixed by the $\hT$-action if and only if
$E$ has a decomposition
\( 
    E = I_1 \oplus \cdots \oplus I_r
\)
satisfying the following conditions for $\alpha = 1,\dots, r$\rom:
\begin{aenume}
\item $I_\alpha$ is an ideal sheaf of $0$-dimensional subscheme
$Z_\alpha$ contained in $\C^2 = \proj^2\setminus\linf$.
\item Under $\Phi$, $I_\alpha|_{\linf}$ is mapped to the $\alpha$-th
factor $\shfO_{\linf}$ of $\shfO_{\linf}^{\oplus r}$.
\item $I_\alpha$ is fixed by the action of $\C^*\times\C^*$, coming from
that on $\proj^2$.
\end{aenume}

\rom{(2)} The fixed point set consists of finitely many points
parametrized by $r$-tuple $(Y_1, \dots, Y_r)$ of Young diagrams such that
\(
   \sum_\alpha |Y_\alpha| = n
\)
\rom(by a way explained in the proof\rom).

\rom{(3)} The fixed point set in $M_0(r,n)$ 
\rom(more strongly, the fixed point set with respect to the first two
factors $T^2$ of $T^{r+2}$\rom)
consists of a single point
$n[0] \in S^n\C^2 \subset M_0(r,n)$.
\end{Proposition}

\begin{proof}
(1) $E\in M(r,n)$ is fixed by the latter $T$-action if and only if
it decomposes as $E = I_1\oplus\dots\oplus I_r$ ($I_\alpha\in
M(1,n_\alpha)$) such that $I_\alpha|_{\linf}$ is mapped to the
$\alpha$-th factor $\shfO_{\linf}$ of $\shfO_{\linf}^{\oplus r}$ under
$\Phi$. Since the double dual $I_\alpha^{\vee\vee}$ of $I_\alpha$ is a
line bundle with $c_1(I_\alpha^{\vee\vee}) = 0$, it is the structure
sheaf $\shfO_{\proj^2}$. Via the natural inclusion $I_\alpha\subset
I_\alpha^{\vee\vee} = \shfO_{\proj^2}$, $I_\alpha$ is an ideal sheaf
of $0$-dimensional subscheme $Z_\alpha$ contained in $\C^2$. Thus
conditions~a),b) are met for $I_\alpha$.
If $E$ is fixed also by the first $T^2$-action, then the condition~c)
must be satisfied. The converse is clear.

(2) By a result of Ellingsrud and Str\o mme~\cite{ES} (see
\cite[\S5.2]{Lecture}) that $I_\alpha$ is fixed if and only if it is
generated by monomials $x^i y^j$, where we consider $I_\alpha$ as an
ideal of $\C[x,y]$, the coordinate ring of $\C^2$. Thus $I_\alpha$
corresponds to a Young diagram $Y_\alpha$ by the rule indicated by the
figure~\ref{fig:young}.
(A monomial $x^{i-1}y^{j-1}$ is placed at $(i,j)$. The ideal
$I_\alpha$ is linearly spanned by monomials outside the Young diagram
$Y_\alpha$. Note that our Young diagrams are rotated $90^\circ$ from
ones used in \cite{Macdonald}.)

\begin{figure}[htbp]
\begin{center}
\leavevmode
\setlength{\unitlength}{0.0005in}
{\renewcommand{\dashlinestretch}{30}
\begin{picture}(3754,3921)(0,-10)
\path(333,3375)(933,3375)(933,2775)
        (333,2775)(333,3375)
\path(333,2775)(933,2775)(933,2175)
        (333,2175)(333,2775)
\path(333,2175)(933,2175)(933,1575)
        (333,1575)(333,2175)
\path(333,1575)(933,1575)(933,975)
        (333,975)(333,1575)
\path(933,1575)(1533,1575)(1533,975)
        (933,975)(933,1575)
\path(933,2175)(1533,2175)(1533,1575)
        (933,1575)(933,2175)
\path(1533,2175)(2133,2175)(2133,1575)
        (1533,1575)(1533,2175)
\path(1533,1575)(2133,1575)(2133,975)
        (1533,975)(1533,1575)
\path(2133,1575)(2733,1575)(2733,975)
        (2133,975)(2133,1575)
\path(333,975)(933,975)(933,375)
        (333,375)(333,975)
\path(933,975)(1533,975)(1533,375)
        (933,375)(933,975)
\path(1533,975)(2133,975)(2133,375)
        (1533,375)(1533,975)
\path(2133,975)(2733,975)(2733,375)
        (2133,375)(2133,975)
\path(2733,975)(3333,975)(3333,375)
        (2733,375)(2733,975)
\thicklines
\blacken\thinlines
\path(63.000,3555.000)(33.000,3675.000)(3.000,3555.000)(63.000,3555.000)
\path(33,3675)(33,75)(3633,75)
\blacken\path(3513.000,45.000)(3633.000,75.000)(3513.000,105.000)(3513.000,45.000)
\put(3708,0){\makebox(0,0)[lb]{\smash{$i$}}}
\put(33,3750){\makebox(0,0)[lb]{\smash{$j$}}}
\put(3558,600){\makebox(0,0)[lb]{\smash{$x^5$}}}
\put(2358,1800){\makebox(0,0)[lb]{\smash{$x^3 y^2$}}}
\put(1158,2400){\makebox(0,0)[lb]{\smash{$x y^3$}}}
\put(558,3600){\makebox(0,0)[lb]{\smash{$y^5$}}}
\put(3033,1200){\makebox(0,0)[lb]{\smash{$x^4 y$}}}
\end{picture}
}
\caption{Young diagram and ideal}
\label{fig:young}
\end{center}
\end{figure}
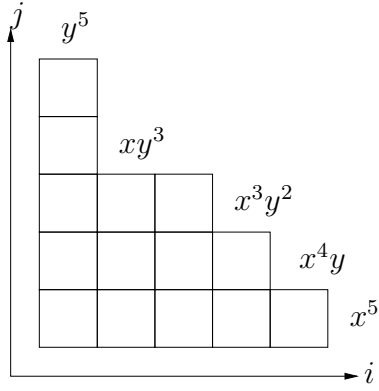

(3) Let us use the description \eqref{eq:alg_quot}. Suppose that the
equivalence class of $(B_1,B_2,i,j)$ is fixed by the $T^2$-action. We
may assume that $(B_1,B_2,i,j)$ has a closed $\GL_n(\C)$-orbit. Then
the equivalence class is fixed if and only if 
\(
   (t_1 B_1, t_2 B_2, i, t_1 t_2 j)
\)
lies in the same $\GL_n(\C)$-orbit.
Since $(t_1 B_1, t_2 B_2, i, t_1 t_2 j)$ converges to
$(0,0,i,0)$ when $t_1, t_2 \to 0$, $(0,0,i,0)$ lies in the closure of
the orbit.
But the orbit is closed, so $(0,0,i,0)$ must be in the orbit.
But $\GL_n(\C)\cdot (0,0,i,0)$ is closed if and only if $i = 0$.
Hence we have $(B_1,B_2,i,j) = (0,0,0,0)$.
\end{proof}

We denote by $\vec{Y}$ an $r$-tuple of Young diagrams $(Y_1,\dots,
Y_r)$. We write the number of boxes of $Y_\alpha$ by $|Y_\alpha|$ and
we set $|\vec{Y}| = \sum_\alpha |Y_\alpha|$.

Let $T_{(E,\Phi)} M(r,n)$ be the tangent space of $M(r,n)$ at a point
$(E,\Phi)$. If $(E,\Phi)$ is fixed by the torus action, then
$T_{(E,\Phi)} M(r,n)$ is a module of the torus. In order to express
the module structure in terms of Young diagrams $Y_\alpha$, we
introduce the following notation.
Let $Y = (\lambda_1\ge \lambda_2 \ge \dots)$ be a Young diagram,
where $\lambda_i$ is the length of the $i$th column.
Let $Y' = (\lambda'_1\ge \lambda_2' \ge \dots)$ be the transpose of
$Y$. Thus $\lambda'_j$ is the length of the $j$th row of $Y$.
Let $l(Y)$ denote the number of columns of $Y$, i.e., $l(Y) =
\lambda'_1$.
Let
\begin{alignat*}{2}
& a_Y(i,j) = \lambda_i - j, & \qquad & a'(i,j) = j - 1 \\
& l_Y(i,j) = \lambda'_j - i, &\qquad & l'(i,j) = i - 1.
\end{alignat*}
Here we set $\lambda_i = 0$ when $i > l(Y)$. Similarly $\lambda'_j =
0$ when $j > l(Y')$.
When the square $s = (i,j)$ lies in $Y$, these are called {\it
arm-length}, {\it arm-colength}, {\it leg-length}, {\it leg-colength}
respectively, and we usually consider in this case. But our formula
below involves these also for squares outside $Y$. So these
take negative values in general. Note that $a'$ and $l'$ does not
depend on the diagram, and we do not write the subscript $Y$.

If two Young diagrams $Y_\alpha$ and $Y_\beta$ are given, we separate
$Y_\alpha$ into two regions ${}^{\heartsuit}Y_\alpha$ and
${}^{\spadesuit}Y_\alpha$ as
\begin{equation*}
  {}^{\heartsuit}Y_\alpha = \left\{ (i,j) \in Y_\alpha \mid j\le
 l(Y_\beta') \right\},
 \qquad
  {}^{\spadesuit}Y_\alpha = \left\{ (i,j) \in Y_\alpha \mid j >
 l(Y_\beta') \right\}.
\end{equation*}
If $l(Y_\alpha') \le l(Y_\beta')$, then ${}^{\spadesuit}Y_\alpha =
\emptyset$. Exchanging the role of $\alpha$ and $\beta$, we divide
$Y_\beta$ into ${}^{\heartsuit}Y_\beta$ and ${}^{\spadesuit}Y_\beta$.
Note that either ${}^{\spadesuit}Y_\alpha$ or ${}^{\spadesuit}Y_\beta$
is the empty set.

\begin{Notation}\label{not:module}
We denote by $e_\alpha$ ($\alpha=1,\dots, r$) the one dimensional
$\hT$-module given by
\begin{equation*}
   \hT\ni (t_1,t_2, e_1, \dots, e_r) \mapsto e_\alpha.
\end{equation*}
Similarly, $t_1$, $t_2$ denote one-dimensional $\hT$-modules. Thus
the representation ring $R(\hT)$ is isomorphic to $\Z[t_1^\pm,
t_2^\pm, e_1^\pm, \dots, e_r^\pm]$, where $e_\alpha^{-1}$ is the dual
of $e_\alpha$.
\end{Notation}

\begin{Theorem}\label{thm:M(r,n)weights}
Let $(E,\Phi)$ be a fixed point of $\hT$-action corresponding to
$\vec{Y} = (Y_1,\dots, Y_r)$. Then the $\hT$-module structure of
$T_{(E,\Phi)} M(r,n)$ is given by
\begin{equation*}
   T_{(E,\Phi)} M(r,n) 
   = \sum_{\alpha,\beta=1}^r N_{\alpha,\beta}(t_1,t_2),
\end{equation*}
where
\[
   N_{\alpha,\beta}(t_1,t_2) = e_\beta\, e_\alpha^{-1}\times
   \left\{
   \sum_{s \in Y_\alpha}
      \left( t_1^{-l_{Y_\beta}(s)} t_2^{a_{Y_\alpha}(s)+1}\right)
%
    + \sum_{t\in Y_\beta} 
        \left(t_1^{l_{Y_\alpha}(t)+1}t_2^{-a_{Y_\beta}(t)}\right)
        \right\}
.
\]
\end{Theorem}

\begin{Remark}
  \textup{(1)} After the first version of this paper was written, the
  authors noticed that this formula already appeared in the context of
  the wall-crossing formula for the Donaldson invariants
  \cite[Lemma~6.2]{EG2}. Their proof does not use the ADHM
  description, so different from ours. We will discuss the relation
  between Nekrasov's prepotential and the wall-crossing formula in a
  future publication with L.~G\"ottsche.
  
  \textup{(2)} The following proof was mentioned also in \cite{BFMT}.
  
\textup{(3)} For the proof of the blowup equation, we only need the
relation between $N_{\alpha,\beta}(t_1,t_2)$ and similar weights on
the blowup (\thmref{thm:bM(r,k,n)weights}). A reader in hurry can
safely skip the proof.
\end{Remark}

\begin{proof}[Proof of \thmref{thm:M(r,n)weights}]
We use the description by linear data for the calculation, which is
very similar to that in \cite[5.8]{Lecture}.

Let $(B_1,B_2,i,j)$ be a datum as above. We consider a complex
\begin{equation}\label{eq:comp} 
   \Hom (V,V) \overset{\sigma}{\longrightarrow}
        \begin{matrix} \Hom(V, V) \oplus \Hom(V, V) \\ 
                        \oplus \\
                       \Hom(W, V) \oplus \Hom(V, W)
                       \end{matrix}
   \overset{\tau}{\longrightarrow} \Hom (V,V),
\end{equation}
where $\sigma$ and $\tau$ are defined by 
\begin{equation*}
   \sigma(\xi) = \begin{pmatrix} \xi B_1 - B_1 \xi \\
                            \xi B_2 - B_2 \xi \\
                            \xi i \\
                            - j\xi \end{pmatrix}, \quad
   \tau \begin{pmatrix} C_1 \\ C_2 \\ I \\ J \end{pmatrix}
        = [B_1,C_2] + [C_1,B_2] + iJ + Ij.
\end{equation*}
This $\sigma$ is the differential of $\GL(V)$-action and $\tau$ is the
differential of the map $(B_1,B_2,i,j)\mapsto [B_1,B_2]+ij$. One can
show that $\sigma$ is injective and $\tau$ is surjective, and the
tangent space of $M(r,n)$ at $\GL(V)\cdot(B_1,B_2,i,j)$ is isomorphic
to the middle cohomology group of the above complex
(cf.\ \cite[1.9]{Lecture} or \cite[3.10]{Quiver}).

Now suppose $\GL(V)\cdot(B_1,B_2,i,j)$ is fixed by the
$\hT$-action. This means that for any $(t_1,t_2,e)\in \hT$
there exists an element $g(t_1,t_2,e)\in \GL(V)$ such that
\begin{equation*}
   (t_1 B_1, t_2 B_2, i e^{-1}, t_1 t_2 e j)
   = g(t_1,t_2,e)^{-1}\cdot (B_1,B_2, i,j).
\end{equation*}
Moreover, such $g(t_1,t_2,e)$ is unique since the $\GL(V)$-action on
the set of stable points is free. In particular, it implies that the
map $(t_1,t_2,e) \mapsto g(t_1,t_2,e)$ is a group homomorphism. We
consider $V$ as a $\hT$-module via it. Also $W$ is a
$\hT$-module via $(t_1,t_2,e)\mapsto e\in \GL(W)$.

We can make the complex~\eqref{eq:comp} {\it $T^{r+2}$-equivariant\/}
by modifying it as
\begin{equation}\label{eq:comp'}
   \Hom (V,V) \overset{\sigma}{\longrightarrow}
        \begin{matrix} t_1 \Hom(V, V) \oplus t_2 \Hom(V, V) \\ 
                        \oplus \\
                        \Hom(W, V) \oplus t_1 t_2 \Hom(V, W)
        \end{matrix}
   \overset{\tau}{\longrightarrow} t_1 t_2 \Hom(V,V),
\end{equation}
where $t_1$, $t_2$ denote the one dimensional $T^{r+2}$-modules as in
Notation~\ref{not:module}.

We have a decomposition
\(
   W = \bigoplus_{\alpha=1}^r e_\alpha
\)
as $T^{r+2}$-modules. From the stability condition, it is easy to see
that $V$ decomposes as $V = \bigoplus V_\alpha e_\alpha$, where $V_\alpha$
is a $T^2$-module (i.e., $T^r$ acts trivially on $V_\alpha$). Thus
$\Ker\tau/\Ima\sigma$ decomposes as
\(
   \bigoplus_{\alpha,\beta}
   \left(\Ker\tau_{\beta\alpha}/\Ima\sigma_{\beta\alpha}\right)
   e_\beta e_\alpha^{-1}
\)
where
\begin{equation}\label{eq:comp''}
   \Hom (V_\alpha,V_\beta) \overset{\sigma_{\beta\alpha}}{\longrightarrow}
        \begin{matrix} t_1 \Hom(V_\alpha, V_\beta) 
                \oplus t_2 \Hom(V_\alpha, V_\beta) \\ 
                        \oplus \\
                           \Hom(W_\alpha, V_\beta)
                \oplus t_1 t_2 \Hom(V_\alpha, W_\beta)
                       \end{matrix}
   \overset{\tau_{\beta\alpha}}{\longrightarrow} 
   t_1 t_2 \Hom(V_\alpha,V_\beta).
\end{equation}

It is clear that each summand has the weight $e_\beta e_\alpha^{-1}$
as a latter torus $T$, so we suppress this factor and only consider
the $\C^*\times\C^*$-module structure hereafter.

Let us write
\(
   Y_\alpha = (\lambda_{\alpha,1}\ge \lambda_{\alpha,2} \ge \dots),
\)
\(
   Y_\alpha' = (\lambda_{\alpha,1}'\ge \lambda_{\alpha,2}'
   \ge \dots).
\)
Since $V_\alpha$ has a basis 
$\{ x^{i-1}y^{j-1} \}$ ($(i,j)\in Y_\alpha$), we have
\begin{equation*}
   V_\alpha = \sum_{j=1}^{\lambda_{\alpha,1}}
             \sum_{i=1}^{\lambda_{\alpha,j}'} t_1^{-i+1} t_2^{-j+1}
           = \sum_{i=1}^{\lambda_{\alpha,1}'}
             \sum_{j=1}^{\lambda_{\alpha,i}} t_1^{-i+1} t_2^{-j+1}.
\end{equation*}
Hence we get
\begin{equation*}
\begin{split}
    & (t_1 + t_2 - 1 - t_1t_2) V_\alpha^* \otimes V_\beta \\
   =\; & \sum_{i=1}^{\lambda_{\alpha,1}'}\sum_{j'=1}^{\lambda_{\alpha,i}}
             t_1^{i-1} t_2^{j'-1} (t_2 - 1) \; \times \;
       \sum_{j=1}^{\lambda_{\beta,1}}\sum_{i'=1}^{\lambda_{\beta,j}'}
      t_1^{-i'+1}(1 - t_1)\, t_2^{-j+1} \\
   =\; & \sum_{i=1}^{\lambda_{\alpha,1}'}\sum_{j=1}^{\lambda_{\beta,1}}
      (t_1^{i-\lambda_{\beta,j}'} - t_1^i)
      (t_2^{-j+\lambda_{\alpha,i}+1} - t_2^{-j+1}) \\
   =\; & \sum_{i=1}^{\lambda_{\alpha,1}'}\sum_{j=1}^{\lambda_{\beta,1}}
      \left[ t_1^{i-\lambda_{\beta,j}'} t_2^{-j+\lambda_{\alpha,i}+1}
       - t_1^i t_2^{-j+1} - (t_1^{i-\lambda_{\beta,j}'} - t_1^i)\,t_2^{-j+1}
       - t_1^i (t_2^{-j+\lambda_{\alpha,i}+1} - t_2^{-j+1})\right].
\end{split}
\end{equation*}
Note that
\begin{equation*}
   \sum_{i=1}^{\lambda_{\alpha,1}'}\sum_{j=1}^{\lambda_{\beta,1}}
      (t_1^{i-\lambda_{\beta,j}'} - t_1^i)\, t_2^{-j+1}
  = \sum_{j=1}^{\lambda_{\beta,1}} \sum_{i=1}^{\lambda_{\beta,j}'}
      (t_1^{1-i} - t_1^{\lambda_{\alpha,1}'-i+1})\, t_2^{-j+1}
  = V_\beta 
     - \sum_{j=1}^{\lambda_{\beta,1}}\sum_{i=1}^{\lambda_{\beta,j}'} 
      t_1^{\lambda_{\alpha,1}'-i+1}t_2^{-j+1}.
\end{equation*}
Similarly note that
\begin{equation*}
   \sum_{i=1}^{\lambda_{\alpha,1}'}\sum_{j=1}^{\lambda_{\beta,1}}
    t_1^i (t_2^{-j+\lambda_{\alpha,i}+1} - t_2^{-j+1})
   = t_1 t_2 V_\alpha^* 
     - \sum_{i=1}^{\lambda_{\alpha,1}'}\sum_{j=1}^{\lambda_{\alpha,i}}
        t_1^i t_2^{-\lambda_{\beta,1}+j} .
\end{equation*}
Thus we have
\begin{equation}\label{eq:N_ba}
\begin{split}
   & \Ker\tau_{\beta\alpha}/\Ima\sigma_{\beta\alpha}
   = (t_1 + t_2 - 1 - t_1t_2) V_\alpha^* \otimes V_\beta
      + V_\beta + t_1 t_2 V_\alpha^* \\
   = \;& \sum_{i=1}^{\lambda_{\alpha,1}'}\sum_{j=1}^{\lambda_{\beta,1}}
     (t_1^{i-\lambda_{\beta,j}'}t_2^{-j+\lambda_{\alpha,i}+1}
        - t_1^i t_2^{-j+1})
    + \sum_{j=1}^{\lambda_{\beta,1}}\sum_{i=1}^{\lambda_{\beta,j}'} 
        t_1^{\lambda_{\alpha,1}'-i+1} t_2^{-j+1}
    + \sum_{i=1}^{\lambda_{\alpha,1}'}\sum_{j=1}^{\lambda_{\alpha,i}}
        t_1^i t_2^{-\lambda_{\beta,1}+j} .
\end{split}
\end{equation}
This is equal to $N_{\alpha,\beta}(t_1,t_2)$, which we want to
compute. We decompose it as
\(
  N_{\alpha,\beta}(t_1,t_2)
  = N^{>0}_{\alpha,\beta}(t_1,t_2) + N^{\le 0}_{\alpha,\beta}(t_1,t_2)
\),
according to the power of $t_2$.
Then
\begin{equation}\label{eq:N+}
\begin{split}
   N^{>0}_{\alpha,\beta}(t_1,t_2) & = 
   \sum_{i=1}^{\lambda_{\alpha,1}'}
   \sum_{j=1}^{\min(\lambda_{\beta,1},\lambda_{\alpha,i})}
     t_1^{i-\lambda_{\beta,j}'} t_2^{-j+\lambda_{\alpha,i}+1}
    + \sum_{i=1}^{\lambda_{\alpha,1}'}
      \sum_{j=\lambda_{\beta,1}+1}^{\lambda_{\alpha,i}}
        t_1^i t_2^{-\lambda_{\beta,1}+j} \\
  & = \sum_{s \in {}^{\heartsuit}Y_\alpha}
       t_1^{-l_{Y_\beta}(s)} t_2^{a_{Y_\alpha}(s)+1} 
  + \sum_{s \in {}^{\spadesuit}Y_\alpha}
       t_1^{l'(s)+1} t_2^{a'(s)-l(Y_\beta')+1} ,
\end{split}
\end{equation}
where the sum $\sum_{j=\lambda_{\beta,1}+1}^{\lambda_{\alpha,i}}$
is understood as $0$ unless $\lambda_{\beta,1}<\lambda_{\alpha,i}$.

Noticing the symmetry 
\(
   N_{\alpha,\beta}(t_1,t_2) = N_{\beta,\alpha}(t_1^{-1},t_2^{-1})t_1t_2
\),
we get the following from \eqref{eq:N_ba}:
\begin{equation*}
   N_{\alpha,\beta}(t_1,t_2)
   = \sum_{i=1}^{\lambda_{\beta,1}'}\sum_{j=1}^{\lambda_{\alpha,1}}
     (t_1^{-i+\lambda_{\alpha,j}'+1} t_2^{j-\lambda_{\beta,i}}
       - t_1^{-i+1}t_2^{j})
    + \sum_{j=1}^{\lambda_{\alpha,1}}\sum_{i=1}^{\lambda_{\alpha,j}'} 
       t_1^{-\lambda_{\alpha,1}'+i} t_2^{j}
    + \sum_{i=1}^{\lambda_{\beta,1}'}\sum_{j=1}^{\lambda_{\beta,i}}
       t_1^{-i+1} t_2^{\lambda_{\alpha,1}-j+1} .
\end{equation*}
This implies that
\begin{equation}\label{eq:N-}
\begin{split}
   N^{\le 0}_{\alpha,\beta}(t_1,t_2) & = 
   \sum_{i=1}^{\lambda_{\beta,1}'}
   \sum_{j=1}^{\min(\lambda_{\alpha,1},\lambda_{\beta,i})}
     t_1^{-i+\lambda_{\alpha,j}'+1} t_2^{j-\lambda_{\beta,i}}
    + \sum_{i=1}^{\lambda_{\beta,1}'}
      \sum_{j=\lambda_{\alpha,1}+1}^{\lambda_{\beta,i}}
        t_1^{-i+1} t_2^{\lambda_{\alpha,1}-j+1}  \\
   &= \sum_{t\in {}^{\heartsuit}Y_\beta} 
        t_1^{l_{Y_\alpha}(t)+1} t_2^{-a_{Y_\beta}(t)}
 + \sum_{t\in {}^{\spadesuit}Y_\beta}
        t_1^{-l'(t)} t_2^{-a'(t)+l(Y_\alpha')} ,
\end{split}
\end{equation}
where the sum $\sum_{j=\lambda_{\alpha,1}+1}^{\lambda_{\beta,i}}$ is
understood as $0$ unless $\lambda_{\alpha,1}<\lambda_{\beta,i}$.
Combining \eqref{eq:N+} with \eqref{eq:N-}, we get the assertion. We
use $-l_{Y_\beta}(s) = l'(s)+1$ for $s\in {}^{\spadesuit}Y_\alpha$ and
re-order the product in $s\in {}^{\spadesuit}Y_\alpha$.
\end{proof}

\section{Moduli spaces on the blowup}\label{sec:mod_blowup}

Let $\bp$ be the blowup of $\proj^2$ at $[1:0:0]$. Let $p\colon
\bp\to\proj^2$ denote the projection. The manifold $\bp$ is the closed
subvariety of $\proj^2\times\proj^1$ defined by
\begin{equation*}
   \{ ([z_0 : z_1 : z_2], [z : w] \in \proj^2\times\proj^1\mid 
   z_1 w = z_2 z \},
\end{equation*}
where the map $p\colon \bp\to\proj^2$ is the projection to the first
factor. Let us denote the inverse image of $\linf$ under
$\bp\to\proj^2$ also by $\linf$ for brevity. It is given by the
equation $z_0 = 0$. The complement $\bp\setminus\linf$ is the blowup
${\widehat\C}^2$ of $\C^2$ at the origin.
Let $C$ denote the exceptional set. It is given by
$z_1 = z_2 = 0$.

In this section, $\shfO$ denotes the structure sheaf of $\bp$,
$\shfO(C)$ the line bundle associated with the divisor $C$,
$\shfO(mC)$ its $m$th tensor power.

Let $\bM(r,k,n)$ be the framed moduli space of torsion free sheaves
$(E,\Phi)$ on $\bp$ with rank $r$, $\langle c_1(E),[C]\rangle = -k$
and $\langle c_2(E) - \frac{r-1}{2r} c_1(E)^2, [\bp]\rangle = n$.

By the same argument as in \propref{prop:vanish} we have
$\Hom(E, E(-\linf)) = \Ext^2(E, E(-\linf)) = 0$ and $\bM(r,k,n)$ is a
nonsingular variety of dimension $2nr$.

\begin{Theorem}\label{thm:blowup}
There is a projective morphism $\widehat\pi\colon\bM(r,k,n)\to
M_0(r,n-\frac1{2r}k(r-k))$ \textup($0\le k<r$\textup) defined by
\begin{equation*}
   (E,\Phi) \mapsto \left(((p_* E)^{\vee\vee}, \Phi),
   \Supp(p_*E^{\vee\vee}/p_*E) + \Supp(R^1p_* E)\right).
\end{equation*}
\end{Theorem}

The proof of this result will be given in \cite{lecture}. In fact, we
prove also the corresponding result for arbitrary projective
surfaces. For the above case with $k = 0$, we can use King's result
\cite{King} instead. Namely there is a morphism from the Uhlenbeck
(partial) compactification $\bM_0(r,0,n) \to M_0(r,n)$ defined via the
ADHM descriptions of both spaces. Then we compose the morphism
$\bM(r,0,n)\to \bM_0(r,0,n)$. This morphism can be defined via a
modification of King's description as in the case of $\C^2$.

We use this result to prove the vanishing result
(\propref{prop:lowblowup}), which is about the case $k=0$, and we can
avoid its usage for the definition of the partition function
$\widehat{Z}$ on the blowup. In this sense, this paper does not rely
on \cite{lecture}.

Let us define an action of the $(r+2)$-dimensional torus $\hT =
\C^*\times\C^* \times T$ on $\bM(r,k,n)$ by modifying the action on
$M(r,n)$ as follows. For $(t_1,t_2)\in \C^*\times\C^*$, let
$F'_{t_1,t_2}$ be an automorphism of $\bp$ defined by
\[
    F'_{t_1,t_2}([z_0: z_1 : z_2], [z : w]) 
    = ([z_0: t_1 z_1 : t_2 z_2], [t_1 z : t_2 w]).
\]
Then we define the action by replacing $F_{t_1,t_2}$ by $F'_{t_1,t_2}$
in \eqref{eq:action}. The action of the latter $T$ is exactly the same 
as before. The morphism $\widehat\pi$ is equivariant.

Note that the fixed point set of $\C^*\times\C^*$ in ${\widehat\C}^2 =
\bp \setminus\linf$ consists of two points $([1 : 0 : 0], [1 : 0])$,
$([1 : 0 : 0], [0 : 1])$. Let us denote them $p_1$ and $p_2$.

Since $C$ is invariant under the $\C^*\times\C^*$-action, the
corresponding line bundle $\shfO(C)$ is an equivariant line bundle.
The section $z_1 / z = z_2 / w$ is equivariant.

\begin{Proposition}\label{prop:fixedpoint2}
\rom{(1)} $(E,\Phi)\in \bM(r,k,n)$ is fixed by the $\hT$-action if
and only if $E$ has a decomposition 
\( 
   E = I_1(k_1 C) \oplus \cdots \oplus I_r(k_r C)
\)
satisfying the following conditions for $\alpha = 1,\dots, r$\rom:
\begin{aenume}
\item $I_\alpha(k_\alpha C)$ is the tensor product
$I_\alpha\otimes\shfO(k_\alpha C)$, where $k_\alpha\in\Z$ and
$I_\alpha$ is an ideal sheaf of $0$-dimensional subscheme $Z_\alpha$
contained in ${\widehat\C}^2 = \bp\setminus\linf$.
\item Under $\Phi$, $I_\alpha(k_\alpha C)|_{\linf}$ is mapped to the
$\alpha$-th factor $\shfO_{\linf}$ of $\shfO_{\linf}^{\oplus r}$.
\item $I_\alpha$ is fixed by the action of $\C^*\times\C^*$, coming
from that on $\bp$.
\end{aenume}

\rom{(2)} The fixed point set consists of finitely many points
parametrized by $r$-tuples 
\(
   ((k_1,Y_1^1, Y_1^2),\linebreak[2]\dots,\linebreak[1] (k_r, Y_r^1, Y_r^2)),
\)
where $k_\alpha\in\Z$ and $Y_\alpha^1$, $Y_\alpha^2$ are
Young diagrams such that
\begin{equation}\label{eq:cond}
   \sum_\alpha k_\alpha = k, \qquad
   \sum_\alpha (|Y_\alpha^1| + |Y_\alpha^2|) + 
   \frac1{2r} \sum_{\alpha < \beta} |k_\alpha - k_\beta|^2 = n
\end{equation}
\rom(by a way explained in the proof\rom).
\end{Proposition}

\begin{proof}
The proof is almost the same as that of \propref{prop:fixedpoint}.

$E\in\bM(r,k,n)$ is fixed by the latter $T^r$-action if and only
if it decomposes as $E = E_1\oplus \dots \oplus E_r$
($E_\alpha\in\bM(1,k_\alpha,n_\alpha)$) such that $E_\alpha|_{\linf}$
is mapped to the $\alpha$-th factor $\shfO_{\linf}$ of
$\shfO_{\linf}^{\oplus r}$ under $\Phi$. Since the double dual
$E_\alpha^{\vee\vee}$ is a line bundle which is trivial at $\linf$, it
is equal to $\shfO(k_\alpha C)$ for some $k_\alpha\in\Z$. Thus
$E_\alpha$ is equal to $I_\alpha(k_\alpha C) = I_\alpha\otimes
\shfO(k_\alpha C)$ for some ideal sheaf $I_\alpha$ of $0$-dimensional
subscheme $Z_\alpha$ in ${\widehat\C}^2$.

If $E$ is fixed also by the first $T^2$-ation, then $I_\alpha$ (and
$Z_\alpha$) is fixed. The support of $Z_\alpha$ must be contained in
the fixed point set in ${\widehat\C}^2$, i.e., $\{ p_1, p_2 \}$.
Thus $Z_\alpha$ is a union of $Z_\alpha^1$ and $Z_\alpha^2$,
subschemes supported at $p_1$ and $p_2$ respectively.
If we take a coordinate system $(x, y) = (z_1/z_0, w/z)$ (resp.\ $=
(z/w, z_2/z_0)$) around $p_1$ (resp.\ $p_2$), then $Z_\alpha^1$
(resp.\ $Z_\alpha^2$) is generated by monomials $x^i y^j$. Then
$Z_\alpha^1$ (resp.\ $Z_\alpha^2$) corresponds to a Young diagram
$Y_\alpha^1$ (resp.\ $Y_\alpha^2$) as before.
\end{proof}

We denote by $\vec{k}$ (resp.\ $\vec{Y}^i$ ($i=1,2$)) for the
$r$-tuple $(k_1,\dots,k_r)$ (resp.\ $(Y^i_1,\dots, Y^i_r$)) as
before. Thus the fixed point corresponds to $(\vec{k}, \vec{Y}^1,
\vec{Y}^2)$.

As in \thmref{thm:M(r,n)weights}, the tangent space of $\bM(r,k,n)$ at
a fixed point $(E,\Phi)$ is a $\hT$-module. 

\begin{Theorem}\label{thm:bM(r,k,n)weights}
Let $(E,\Phi)$ be a fixed point of $\hT$-action corresponding to
$(\vec{k}, \vec{Y}^1,\vec{Y}^2)$. Then the $\hT$-module structure of
$T_{(E,\Phi)} \bM(r,k,n)$ is given by
\begin{equation*}
   T_{(E,\Phi)} \bM(r,k,n)
   = \sum_{\alpha,\beta=1}^r \left(L_{\alpha,\beta}(t_1,t_2) 
     + t_1^{k_\beta - k_\alpha} M^1_{\alpha,\beta}(t_1,t_2)
     + t_2^{k_\beta - k_\alpha} M^2_{\alpha,\beta}(t_1,t_2)\right),
\end{equation*}
where
\begin{equation*}
   L_{\alpha,\beta}(t_1,t_2) = 
   e_\beta\, e_\alpha^{-1}\times
   \begin{cases}
     {\displaystyle
     \sum_{\substack{i,j\ge 0\\i+j \le k_\alpha-k_\beta-1}}
          t_1^{-i} t_2^{-j}}
       & \text{if $k_\alpha > k_\beta$}, \\
     {\displaystyle
     \sum_{\substack{i,j\ge 0\\i+j \le k_\beta-k_\alpha-2}}
          t_1^{i+1} t_2^{j+1}}
       & \text{if $k_\alpha + 1 < k_\beta$}, \\
     0 & \text{otherwise},
   \end{cases}
\end{equation*}
and $M^1_{\alpha,\beta}(t_1,t_2)$
\rom(resp.\ $M^2_{\alpha,\beta}(t_1,t_2)$\rom)
is equal to $N_{\alpha,\beta}(t_1,t_2/t_1)$
\rom(resp.\ $N_{\alpha,\beta}(t_1/t_2, t_2)$\rom),
with $(Y_\alpha, Y_\beta)$ is replaced by $(Y_\alpha^1, Y_\beta^1)$
\rom(resp.\ $(Y_\alpha^2, Y_\beta^2)$\rom).
\end{Theorem}

\begin{proof}
According to the decomposition $E = I_1(k_1 C)\oplus\cdots\oplus
I_r(k_r C)$, the tangent space $T_{(E,\Phi)} \bM(r,k,n) = \Ext^1(E,
E(-\linf))$ is decomposed as
\begin{equation*}
   \Ext^1(E, E(-\linf)) 
   = \bigoplus_{\alpha,\beta}
        \Ext^1(I_\alpha(k_\alpha C), I_\beta(k_\beta C - \linf)).
\end{equation*}
The factor $\Ext^1(I_\alpha(k_\alpha C), I_\beta(k_\beta C - \linf))$
has weight $e_\beta e_\alpha^{-1}$ as a $T$-module. Thus our
remaining task is to describe each factor as a $T^2$-module. We
suppress $e_\beta e_\alpha^{-1}$ hereafter.

Let $\Ext^*$ denotes the alternating sum
$\sum_i (-1)^i \Ext^i$ considered as an element of the representation
ring. Then $\Ext^*$ defines a homomorphism from the equivariant
$K$-group to the representation ring.
By \propref{prop:vanish} we have $\Ext^*(I_\alpha(k_\alpha C),
I_\beta(k_\beta C - \linf)) = - \Ext^1(I_\alpha(k_\alpha C),
\linebreak[1] I_\beta(k_\beta C - \linf))$.
Using the exact sequence $0 \to I_\alpha \to \shfO \to
\shfO_{Z_\alpha} \to 0$, we have
\begin{equation}\label{eq:exts}
\begin{split}
   & \Ext^*(I_\alpha(k_\alpha C), I_\beta(k_\beta C - \linf)) \\
 =\;& \Ext^*(\shfO(k_\alpha C), \shfO(k_\beta C - \linf))
  -   \Ext^*(\shfO(k_\alpha C), \shfO_{Z_\beta}(k_\beta C - \linf))\\
  & - \Ext^*(\shfO_{Z_\alpha}(k_\alpha C), \shfO(k_\beta C - \linf))
  +   \Ext^*(\shfO_{Z_\alpha} (k_\alpha C),
           \shfO_{Z_\beta}(k_\beta C - \linf)).
\end{split}
\end{equation}

Let us first consider the term
\(
   \Ext^*(\shfO(k_\alpha C), \shfO(k_\beta C - \linf))
   \linebreak[1]
   = - \Ext^1(\shfO(k_\alpha C), \shfO(k_\beta C - \linf))
   \linebreak[1]
   = - H^1(\shfO((k_\beta - k_\alpha)C - \linf)).
\)
We show that this is equal to $-L_{\alpha,\beta}$.

Set $n = k_\alpha - k_\beta$.
If $n = 0$, we have $H^1(\bp, \shfO(-\linf)) = 0$ by
\propref{prop:vanish}. Thus we have the expression $L_{\alpha,\beta}$
in this case.

Next suppose $n > 0$.
Consider the cohomology long exact sequence associated
with an exact sequence 
\(
   0 \to \shfO(-nC) \to \shfO((-n+1)C) \to \shfO_C((-n+1)C) \to 0
\).
Note that this is equivariant under the $\C^*\times\C^*$-action.
Since $C$ is a projective line $\proj^1$ with self-intersection
$(-1)$, we have 
\( 
   H^1(C, \shfO_C((-n+1)C)) = H^1(\proj^1, \shfO_{\proj^1}(n-1))
   = 0
\).
Thus we have
\begin{equation*}
    0 \to H^0(\proj^1, \shfO_{\proj^1}(n-1)) \to
          H^1(\bp, \shfO(-nC - \linf)) \to
          H^1(\bp, \shfO((-n+1)C - \linf)) \to 0.
\end{equation*}
This is an exact sequence in $\C^*\times\C^*$-modules.
Starting with $H^1(\bp, \shfO(-\linf)) = 0$, we get
\begin{equation*}
    H^1(\bp, \shfO(-nC - \linf)) 
    = \bigoplus_{d=0}^{n-1} H^0(\proj^1, \shfO_{\proj^1}(d))
\end{equation*}
by induction.
Since $H^0(\proj^1, \shfO_{\proj^1}(d))$ is the space of homogeneous
polynomials in $z$, $w$ of degree $d$, it is equal to $\sum_{i=0}^d
t_1^{-i} t_2^{-d+i}$ in the representation ring of $T^2$. Thus we have
the expression $L_{\alpha,\beta}$ in this case.

Finally consider the case $n < 0$. The proof is almost the same as
that for the case $n > 0$. We use
\(
   0 \to \shfO((-n-1)C) \to \shfO(-nC) \to \shfO_C(-nC) \to 0
\)
to get
\begin{equation*}
    0 \to H^1(\bp, \shfO((-n-1)C - \linf)) \to
          H^1(\bp, \shfO(-nC - \linf)) \to
          H^1(\proj^1, \shfO_{\proj^1}(n)) \to 0,
\end{equation*}
where we have used $H^0(\proj^1, \shfO_{\proj^1}(n)) = 0$. Starting
with $H^1(\bp, \shfO((-n-1)C - \linf)) = 0$ for $n=-1$, we get
\(
    H^1(\bp, \shfO(-nC - \linf)) 
    = \bigoplus_{d=1}^{-n} H^1(\proj^1, \shfO_{\proj^1}(-d))
\)
by induction. The canonical bundle $K_{\proj^1}$ of $\proj^1$ is
isomorphic to $\shfO_{\proj^1}(-2)$. But this isomorphism is not
equivariant, and the actual formula is $K_{\proj^1} \cong
t_1^{-1}t_2^{-1}\shfO_{\proj^1}(-2)$.
Therefore the Serre duality says 
\(
   H^1(\proj^1, \shfO_{\proj^1}(-d))
\)
is the dual of
\(
   t_1^{-1}t_2^{-1} H^0(\proj^1, \shfO_{\proj^1}(d-2)).
\)
Thus we get the assertion also in this case.

Now we turn to the remaining three terms in \eqref{eq:exts}. We have
\begin{equation}\label{eq:exts2}
\begin{split}
   &
   \begin{gathered}
     - \Ext^*(\shfO(k_\alpha C), \shfO_{Z_\beta}(k_\beta C - \linf))
    - \Ext^*(\shfO_{Z_\alpha}(k_\alpha C), \shfO(k_\beta C - \linf)) \\
    + \Ext^*(\shfO_{Z_\alpha} (k_\alpha C), \shfO_{Z_\beta}(k_\beta C -
                                \linf))
   \end{gathered}
\\
  =\; & 
  \begin{gathered}[t]
   - \Ext^*(\shfO, \shfO_{Z_\beta}((k_\beta - k_\alpha)C - \linf))
   - \Ext^*(\shfO_{Z_\alpha}((k_\alpha - k_\beta) C), \shfO(- \linf)) \\
   + \Ext^*(\shfO_{Z_\alpha},
              \shfO_{Z_\beta}((k_\beta - k_\alpha)C - \linf)).
  \end{gathered}
\end{split}
\end{equation}

As in the proof of \propref{prop:fixedpoint2}, we have decomposition
$Z_\alpha = Z_\alpha^1 \cup Z_\alpha^2$ according to the support
$p_1$, $p_2$. Hence each of the remaining terms in \eqref{eq:exts2} is
the direct sum of the corresponding terms for $Z_\alpha^1$,
$Z_\beta^1$ and $Z_\alpha^2$, $Z_\beta^2$.
(A mixed term 
\(
   \Ext^*(\shfO_{Z_\alpha^1}(k_\alpha C),
   \shfO_{Z_\beta^2}(k_\beta C - \linf))
\)
is obviously zero.)
We study each summand separately.

First consider terms for $Z_\alpha^1$, $Z_\beta^1$.  We take a
coordinate system $(x, y) = (z_1/z_0, w/z)$ as in the proof of
\propref{prop:fixedpoint2}. Since the divisor $C$ is given by $x = 0$,
the multiplication by $x^m$ induces an isomorphism
$\shfO_{Z_\alpha^1}(mC) \cong \shfO_{Z_\alpha^1}$ of sheaves for
$m\in\Z$. It becomes an isomorphism of {\it equivariant\/} sheaves if
we twist it as $\shfO_{Z_\alpha^1}(mC) \cong t_1^m\shfO_{Z_\alpha^1}$.
Hence the summand of \eqref{eq:exts2} for $p_1$ is equal to
\begin{equation}\label{eq:exts3}
    t_1^{k_\beta - k_\alpha}
    \left(- \Ext^*(\shfO, \shfO_{Z_\beta^1}(- \linf))
          - \Ext^*(\shfO_{Z_\alpha^1}, \shfO(- \linf))
          + \Ext^*(\shfO_{Z_\alpha^1},
                       \shfO_{Z_\beta^1}(- \linf))\right).
\end{equation}
Since $Z_\alpha^1$ is supported at the single point $p_1$, we can
consider it as a subscheme of $\proj^2$ supported at the origin $[1:
0: 0]$, where the $T^2$-action on $\proj^2$ is given by 
\(
  [z_0 : z_1 : z_2]\mapsto [z_0 : t_1 z_1 : t_2/t_1 z_2]
\).
Let $I_\alpha^1$ be the corresponding ideal sheaves of $\shfO_{\proj^2}$.
Using the
\(
  0 \to I_\alpha^1 \to \shfO_{\proj^2} \to \shfO_{Z_\alpha^1}\to 0,
\)
we find that \eqref{eq:exts3} is equal to
\begin{equation*}
    t_1^{k_\beta - k_\alpha}\left(
     \Ext^*(I_\alpha^1, I_\beta^1) 
      - \Ext^*(\shfO_{\proj^2}, \shfO_{\proj^2}(-\linf))
    \right).
\end{equation*}
The second term $\Ext^*(\shfO_{\proj^2}, \shfO_{\proj^2}(-\linf))$
is zero.
Thus we can use the formula in \thmref{thm:M(r,n)weights} after
replacing $(t_1, t_2)$ by $(t_1, t_2/t_1)$, and get the expression
$M^1_{\alpha,\beta}(t_1,t_2)$ in the assertion.

The terms for $Z_\alpha^2$, $Z_\beta^2$ can be calculated in a similar
way. We get $M^2_{\alpha,\beta}(t_1,t_2)$.
\end{proof}

For a future reference, we record the formula of the character:
\begin{Proposition}
\begin{equation*}
   \ch H^1(\bp, \shfO(-kC-\linf))
   = 
   \begin{cases}
     {\displaystyle
     \sum_{\substack{i,j\ge 0\\i+j \le k-1}}
          t_1^{-i} t_2^{-j}}
       & \text{if $k > 0$}, \\
     {\displaystyle
     \sum_{\substack{i,j\ge 0\\i+j \le -k-2}}
          t_1^{i+1} t_2^{j+1}}
       & \text{if $k < -1$}, \\
     0 & \text{if $k=0$ or $-1$}.
   \end{cases}
\end{equation*}
\end{Proposition}

\section{Sums over Young tableaux and Hilbert series}\label{sec:Hilbert}

Although our main concern is about equivariant homology groups of
moduli spaces, equivariant $K$-groups are more natural for an
explanation of meaning of Nekrasov's partition function.

Let $K^{\hT}(M(r,n))$ be the Grothendieck group of $\hT =
T^{r+2}$-equivariant coherent sheaves on $M(r,n)$ and similarly for
$K^{\hT}(\bM(r,k,n))$, $K^{\hT}(M_0(r,n))$. These are modules over the
representation ring $R(\hT)$ of the torus ${\hT}$. As in
\ref{not:module}, we identify it with the Laurent polynomial ring
$\Z[t_1^\pm,t_2^\pm, e_1^\pm\dots, e_r^\pm]$. Since $M(r,n)$ and
$\bM(r,k,n)$ are nonsingular, $K^{\hT}(M(r,n))$, $K^{\hT}(\bM(r,k,n))$
are isomorphic to the Grothendieck groups of ${\hT}$-equivariant
locally free sheaves. In particular, they have the ring structures
given by tensor products. For an equivariant proper morphism $f$
between ${\hT}$-varieties, we have induced homomorphism $f_*$ between
the Grothendieck groups given by the alternating sum of higher direct
image sheaves $\sum_i (-1)^i R^i f_*$. In particular, we have
\begin{equation*}
   \pi_* \colon K^{\hT}(M(r,n))\to K^{\hT}(M_0(r,n)), \qquad
   \widehat\pi_* \colon K^{\hT}(\bM(r,k,n))\to K^{\hT}(M_0(r,n)).
\end{equation*}

Let $\mathcal R = \Q(t_1,t_2,e_1,\dots,e_m)$ be the quotient field of
$R({\hT})$. Let $\iota\colon M(r,n)^{\hT}\to M(r,n)$ be the inclusion
of the ${\hT}$-fixed point set. By the localization theorem for the
$K$-theory due to Thomason~\cite{Th} (a prototype of the localization
theorem was in \cite{AS}), it is known that the homomorphism $\iota_*$
is an isomorphism after the localization:
\begin{equation*}
   \iota_*\colon
   K^{\hT}(M(r,n)^{\hT})\otimes_{R({\hT})}\mathcal R
   \xrightarrow{\cong}
   K^{\hT}(M(r,n))\otimes_{R({\hT})}\mathcal R.
\end{equation*}
Since $M(r,n)^{\hT}$ consists of finitely many points $\{ \vec{Y} \}$, and
$K^{\hT}$ of the point is isomorphic to the representation ring, the left
hand side is the direct sum $\# \{\vec{Y}\}$-copies of $\mathcal
R$. Similarly, $K^{\hT}(\bM(r,k,n))\otimes_{R({\hT})}\mathcal R$ is isomorphic
to $\# \{ (\vec{k}, \vec{Y}^1, \vec{Y}^2)\}$-copies of $\mathcal
R$. On the other hand, $M_0(r,n)^{\hT}$ consists of a single point $\{ 0
\}$, hence $K^{\hT}(M_0(r,n))\otimes_{R({\hT})}\mathcal R \cong\mathcal R$.

The inverse of $\iota_*$ can be explicitly given by the following
formula:
\begin{equation*}
   \iota_*^{-1}(\bullet) =
   \bigoplus_{\vec{Y}} \frac{\iota_{\vec{Y}}^*(\bullet)}
   {\Wedge_{-1}T^*_{\vec{Y}}M(r,n)},
\end{equation*}
where $T^*_{\vec{Y}}M(r,n)$ is the cotangent bundle of $M(r,n)$ at a
fixed point of $\vec{Y}$ considered as a ${\hT}$-module, $\Wedge_{-1}$ is
the alternating sum of exterior powers, and $\iota_{\vec{Y}}^*$ is the
pull-back homomorphism with respect to the inclusion
$\iota_{\vec{Y}}\colon \{ \vec{Y}\}\to M(r,n)$. Here the pull-back
homomorphism is defined via the isomorphism of $K^{\hT}(M(r,n))$ and the
Grothendieck group of ${\hT}$-equivariant {\it locally free sheaves\/}.

If $M(r,n)$ would be compact, we have a pushforward homomorphism
$p_*\colon K^{\hT}(M(r,n))\to R({\hT})$ given by $p\colon M(r,n)\to \{pt\}$
and it can be computed by the Bott's formula:
\begin{equation*}
   p_*(\bullet) = \sum_{\vec{Y}} \frac{\iota_{\vec{Y}}^*(\bullet)}
   {\Wedge_{-1}T^*_{\vec{Y}}M(r,n)}.
\end{equation*}
However $M(r,n)$ is noncompact, and $p_*$ is not defined. In fact,
cohomology groups $H^i(M(r,n), \bullet)$ may be infinite
dimensional. But the right hand side makes sense as an element in
$\mathcal R$. In fact, it computes the alternating sum of {\it Hilbert
series\/} of cohomology groups:
\begin{Proposition}\label{prop:Lef}
Let $E$ be a ${\hT}$-equivariant coherent sheaf on $M(r,n)$. Then we have
\begin{equation*}
   \sum_{i=0}^{2nr} (-1)^i \ch H^i(M(r,n), E)
   = \sum_{\vec{Y}} \frac{\iota_{\vec{Y}}^* E}
   {\Wedge_{-1}T_{\vec{Y}}^*M(r,n)},
\end{equation*}
where $\ch$ denotes the Hilbert series.
\end{Proposition}

Let us recall the definition of the Hilbert series. (See
\cite[\S6.6]{CG} and the reference therein for more detailed account.)
Let $V$ be a representation of the torus ${\hT}$. Let $V = \bigoplus
V_\mu$ ($\mu\in X^*({\hT})$) be its weight space decomposition, i.e.,
\begin{equation*}
   V_\mu = \{ v\in V \mid \text{$t\cdot v = \mu(t) v$ for $t\in {\hT}$} \}.
\end{equation*}
Here $X^*(\hT)$ denotes the group of characters of $\hT$.
When the dimensions of weight spaces are finite dimensional, we define
the character of $V$ by
\begin{equation*}
   \ch V = \sum (\dim V_\mu) e^\mu.
\end{equation*}
We take coordinates $(t_1,t_2,e_1,\dots,e_r)\in {\hT}$ as before, and we
consider the right hand side as an element in the Laurent power series
in $t_1$, $t_2$, $e_1$, \dots, $e_r$.

We want to apply this definition to the cohomology groups of a
${\hT}$-equivariant coherent sheaf on $M(r,n)$. Since $M(r,n)$ is {\it
not\/} projective and cohomology groups are {\it not\/}
finite-dimensional in general, we first need to show that weight
spaces are finite-dimensional and the above definition makes sense.
For this purpose, we consider a ${\hT}$-equivariant coherent sheaf $E$
on the {\it affine algebraic variety\/} $M_0(r,n)$. Then the space $M
= H^0(M_0(r,n), E)$ of global sections of $E$ is identified with a
finitely generated module over the coordinate ring of $M_0(r,n)$. (And
the higher cohomology groups vanishes.) Let $M = \bigoplus M_\mu$
($\mu\in X^*({\hT})$) be the weight space decomposition as above.

\begin{Lemma}
A weight space $M_\mu$ is finite-dimensional as a vector space over
$\C$.
\end{Lemma}

\begin{proof}
By \cite{Lu-qv}, the coordinate ring is generated by the following two
types of elements
\begin{enumerate}
\item $\tr(B_{\alpha_N}B_{\alpha_{N-1}}\cdots B_{\alpha_1}\colon
V \to V)$,
\item $\langle \chi, j B_{\alpha_N} B_{\alpha_{N-1}} \cdots
B_{\alpha_1} i\rangle$,
\end{enumerate}
where $\alpha_1$, \dots, $\alpha_N$ is $1$ or $2$ and $\chi$ is a
linear form on $\End(W)$. Any of these elements is contained in a
weight space with a {\it nonzero\/} weight. From this we get our
assertion.
\end{proof}

Now the Hilbert series of $E$ (or $M$) is defined by
\begin{equation*}
    \ch E \equiv \ch M = \sum_\mu (\dim M_\mu) e^\mu.
\end{equation*}
By a well-known argument on Hilbert series, one can show that $\ch E$
is a rational function, i.e., an element in $\mathcal R$.

Now we can return to the situation in \propref{prop:Lef}. Let $E$ be a
${\hT}$-equivariant coherent sheaf on $M(r,n)$. Since $\pi\colon M(r,n)\to
M_0(r,n)$ is a projective morphism, the higher direct image sheaves
$R^i\pi_* E$ is an equivariant coherent sheaf on $M_0(r,n)$. The space
of its global sections is the higher cohomology group $H^i(M(r,n),E)$.
Thus we can consider the associated Hilbert series
\begin{equation*}
   \ch R^i\pi_* E \equiv \ch H^i(M(r,n),E).
\end{equation*}

Now we can finish the proof of \propref{prop:Lef} thanks to a general
result of Thomason~\cite{Th}. The argument appears in \cite{Haiman}
for $r=1$, and his argument can be applied to our situation, once the
above property of the coordinate ring of $M_0(r,n)$ is established.

The proof follows from the commutativity of the following square
\begin{equation*}
\begin{CD}
   K^{\hT}(M(r,n))\otimes_{R({\hT})}\mathcal R
   @>\cong>(\iota_*)^{-1}>\bigoplus_{\vec{Y}} \mathcal R
\\
   @V\pi_*VV @VV\sum_{\vec{Y}}V
\\
   K^{\hT}(M_0(r,n))\otimes_{R({\hT})}\mathcal R
   @>\cong>(\iota_{0*})^{-1}>\mathcal R
\end{CD}
\end{equation*}
and the observation $(\iota_{0*})^{-1} = \ch$, which is a consequence
of a trivial identity $\ch\circ\iota_{0*} = \operatorname{id}$. Here
$\iota_0$ is the inclusion of the unique fixed point of $M_0(r,n)$.

Let us give two examples. Let $\shfO$ be the structure sheaf of
$M(2,1)$. We directly check that \propref{prop:Lef} holds for $E =
\shfO$.
We have two fixed points $\vec{Y} = ([1],[\emptyset])$,
$([\emptyset],[1])$ in $M(2,1)$. The localization gives us
\begin{equation}\label{eq:M21}
\begin{split}
   & \frac1{(1-t_1)(1-t_2)(1-\frac{e_1}{e_2})(1-t_1 t_2\frac{e_2}{e_1})}
   + 
   \frac1{(1-t_1)(1-t_2)(1-\frac{e_2}{e_1})(1-t_1 t_2\frac{e_1}{e_2})}
\\
  =\; &
  \frac{1 + t_1 t_2}{(1-t_1)(1-t_2)
  (1- t_1t_2\frac{e_1}{e_2})(1- t_1t_2\frac{e_2}{e_1})}.
\end{split}
\end{equation}
On the other hand, we have $M(2,1) \cong \C^2\times T^*\proj^1$. The
$\C^2$-component is given by $(B_1,B_2)$ and $T^*\proj^1$-component is
given by $(\Ker i, j i)$, where $\Ker i$ is a one-dimensional subspace
in the two-dimensional space $W$, and $\xi = ji$ is an endomorphism of
$W$ satisfying $\xi(\Ker i) = 0$, $\Ima\xi\subset\Ker i$. The higher
cohomology groups $H^i(M(2,1), \shfO) = 0$ ($i > 0$) vanish, and the
global sections $H^0(M(2,1),\shfO)$ is identified with
\begin{equation*}
   \C[x,y]\otimes (\C[s,t,u]/st = u^2),
\end{equation*}
where $x = B_1$, $y = B_2$, $s = j_1 i_2$, $t = j_2 i_1$, $u = j_1 i_1 
= - j_2 i_2$ with 
\(
  i = \left[
  \begin{smallmatrix}
    i_1 & i_2
  \end{smallmatrix}
  \right]
\),
\(
  j = \left[
  \begin{smallmatrix}
    j_1 \\ j_2
  \end{smallmatrix}
  \right]
\).
Since weights of $x$, $y$, $s$, $t$, $u$ are $t_1$, $t_2$,
$t_1t_2e_1/e_2$, $t_1t_2e_2/e_1$, $t_1t_2$ respectively, the character 
of $H^0(M(2,1),\shfO)$ is also given by \eqref{eq:M21}.

\begin{Remark}
  We have used the following convention on the action on the
  coordinate ring. Let $F_g\colon M(r,n)\to M(r,n)$ be the isomorphism
  given an element $g\in \hT$. It induces a map $F_g^*$ given by
  $f\mapsto f\circ F_g$ on the coordinate ring. The same applies to
  $H^i(M(r,n), E)$ for a $\hT$-equivariant sheaf $E$.
  Accordingly when we apply \propref{prop:Lef}, we make $\hT$ acts on
  the cotangent space $T^*_{\vec{Y}}M(r,n)$ by $d(F_g)_{\vec{Y}}^*$.
\end{Remark}

Next consider the rank $1$ case. The moduli space $M(1,n)$ is nothing
but the Hilbert scheme $(\C^2)^{[n]}$ of $n$ points in $\C^2$. We
apply \propref{prop:Lef} to the structure sheaf $\shfO$ of $M(1,n)$.
The fixed points are parametrized by Young diagrams $Y$ of size $n$ as
\propref{prop:fixedpoint}. The weights of tangent spaces at fixed
points is given by the formula \ref{thm:M(r,n)weights}. In particular,
the localization gives us
\begin{equation*}
   \sum_{|Y|=n}
   \frac1{\prod_{s\in Y}
   (1 - t_1^{-l(s)} t_2^{1+a(s)})(1-t_1^{1+l(s)}t_2^{-a(s)})}.
\end{equation*}
On the other hand, we have $H^0((\C^2)^{[n]},\shfO) = H^0(S^n(\C^2),
\shfO) = \C[x_1,y_1,\dots, x_n,y_n]^{S_n}$, where $S_n$ acts by
permuting $(x_1,y_1)$, \dots, $(x_n, y_n)$. Higher cohomology groups
$H^i((\C^2)^{[n]}, \shfO)$ ($i > 0$) vanish since $S^n\C^2$ is a
rational singularity. Now $\C[x_1,y_1,\dots, x_n,y_n]^{S_n}$ is
isomorphic to $S^n(\C[x,y])$, and the generating function of the
Hilbert series is given by
\begin{equation*}
\begin{split}
    & \sum_{n=0}^\infty \q^n\, \ch H^0(S^n(\C^2), \shfO)
   = \prod_{p_1,p_2\ge 0} \frac1{1 - t_1^{p_1}t_2^{p_2}\q}
   = \exp \left( - \sum_{p_1,p_2\ge 0} \log(1 - t_1^{p_1}t_2^{p_2}\q)\right)
\\
   =\; & \exp\left(\sum_{p_1,p_2\ge 0}\sum_{r=1}^\infty
     \frac{t_1^{rp_1}t_2^{rp_2}\q^r}r\right)
   = \exp\left(\sum_{r=1}^\infty \frac{\q^r}{(1-t_1^r)(1-t_2^r)r}\right).
\end{split}
\end{equation*}
Thus we get
\begin{equation}\label{eq:Haiman}
   \sum_{Y}
   \frac{\q^{|Y|}}{\displaystyle\prod_{s\in Y}
   (1 - t_1^{-l(s)} t_2^{1+a(s)})(1-t_1^{1+l(s)}t_2^{-a(s)})}
   = \exp\left(\sum_{r=1}^\infty \frac{\q^r}{(1-t_1^r)(1-t_2^r)r}\right).
\end{equation}
A purely combinatorial proof of this identity can be found in
\cite[VI]{Macdonald}.
A different geometric proof can be found in
\cite[Lemma~3.2]{Haiman}. It also uses geometry of Hilbert schemes.

Let $H^{\hT}_*(M(r,n))$ be the $\hT$-equivariant Borel-Moore homology
group of $M(r,n)$ with rational coefficients. We define it as in
\cite[\S2.8]{Lu:GH}, but we assign the degree as in \cite{EG} so that the
fundamental class $[M(r,n)]$ has degree $2\dim M(r,n) = 4rn$.

Let us recall the definition briefly.
We have a finite dimensional approximation of the classifying space
$E\hT\to B\hT$, i.e., for any $n$, there exists a smooth irreducible
variety $U$ with $\hT$-action such that
\begin{aenume}
\item The quotient $U \to U/\hT$ exists and is a principal $\hT$-bundle.
\item $H^i(U) = 0$ for $i=1,\dots,n$.
\end{aenume}
We then define
\begin{equation*}
   H^{\hT}_n(X) = H_{n-2\dim \hT+2\dim U}(X\times_{\hT} U),
\end{equation*}
where $H_*(\ )$ in the right hand side is the Borel-Moore homology
group (see e.g., \cite[\S B.2]{Fulton}).
Note that $U$ is smooth, and $\dim U$ makes sense. One can show
that this is independent of the choice of $U$, using the double
fibration argument.

The equivariant homology group is a module over the usual equivariant
cohomology of a point $H^*_{\hT}(pt)$. The latter is the symmetric
algebra of the dual of the Lie algebra of $\hT$, which we denote by
$S(\hT)$. We choose its generators $\ve_1$, $\ve_2$, $a_1$, \dots,
$a_r$ corresponding to $t_1$, $t_2$, $e_1$, \dots, $e_r$ respectively.
We use the vector notation $\vec{a}$ for $(a_1,\dots,a_r)$.  We have
$H^{\hT}_k(X) = 0$ if $k > 2\dim X$, but $H^{\hT}_k(X)$ may be nonzero
for $k<0$.

The results given in this section have counterparts for equivariant
homology groups. For example, we have a commutative diagram
\begin{equation*}
\begin{CD}
   H^{\hT}_*(M(r,n))\otimes_{S(\hT)}\mathcal S
   @>\cong>(\iota_*)^{-1}>\bigoplus_{\vec{Y}} \mathcal S
\\
   @V\pi_*VV @VV\sum_{\vec{Y}\bullet}V
\\
   H^{\hT}_*(M_0(r,n))\otimes_{S({\hT})}\mathcal S
   @>\cong>(\iota_{0*})^{-1}>\mathcal S,
\end{CD}
\end{equation*}
where $\mathcal S$ is the quotient field of $S(\hT)$. The proof of the
localization theorem for equivariant Borel-Moore homology can be
found, for example, in \cite[4.4]{Lu:GH2}.
We have
\begin{equation*}
   \sum_{\vec{Y}} \frac{\iota_{\vec{Y}}^*\alpha}{e(T_{\vec{Y}})}
   = (\iota_{0*})^{-1}\pi_*(\alpha).
\end{equation*}

Further more, the right hand side has an interpretation as the
{\it equivariant Hilbert polynomial\/} of $E$ if $\alpha$ is the Chern 
character of a vector bundle $E$. For example, we have the following
for $E = \shfO$:
\begin{equation}\label{eq:HilbPoly}
\begin{split}
   & (\iota_{0*})^{-1} \pi_* [M(r,n)]
   = \sum_{\vec{Y}} \frac{1}{e(T_{\vec{Y}})}
   = \lim_{t\to 0}\left.
   \sum_{\vec{Y}} \frac{t^{2nr}}{\Wedge_{-1}T_{\vec{Y}}^*}
 \right|_{\substack{t_1 = e^{-t\ve_1}, t_2 = e^{-t\ve_2}\\
 e_\alpha = e^{-ta_\alpha}}}
\\
   =\; &\lim_{t\to 0} \left. t^{2nr}
   \sum_{i=0}^{2nr} (-1)^i \ch H^i(M(r,n), \shfO)
 \right|_{\substack{t_1 = e^{-t\ve_1}, t_2 = e^{-t\ve_2}\\
 e_\alpha = e^{-ta_\alpha}}}.
\end{split}
\end{equation}
This is our interpretation of Nekrasov's partition function mentioned
in the introduction.
For example, for $r = 1$, we can derive the following from
\eqref{eq:Haiman}:
\begin{equation}\label{eq:formula1}
   \sum_{Y}
   \frac{\q^{|Y|}}{\displaystyle\prod_{s\in Y}
   \left\{ -l_Y(s)\ve_1 + (1+a_Y(s))\ve_2\right\}
   \left\{ (1+l_Y(s))\ve_1 - a_Y(s)\ve_2\right\}}
   = \exp\left(\frac{\q}{\ve_1\ve_2}\right).
\end{equation}
As in the proof of \eqref{eq:Haiman}, we can directly obtain the right
hand side as follows. We use localization on $M_0(1,n) = S^n(\C^2)$,
instead of $M(1,n) = (\C^2)^{[n]}$. The point is that $S^n(\C^2)$ is
an orbifold, and hence has an explicit formula of $(\iota_{0*})^{-1}$.
This formula justifies the following definition of `generating
spaces':
\begin{equation*}
   \exp(\q \C^2) = \sum_{n=0}^\infty \q^n S^n(\C^2),
   \qquad\text{or }
   \q \C^2 = \log\left(\sum_{n=0}^\infty \q^n S^n(\C^2)\right).
\end{equation*}

\section{Rank $1$ case}\label{sec:rank1}

This section is a detour. We study Nekrasov's partition function and
its analog for blowup in the rank $1$ case.

The partition function for rank $1$ is
\begin{equation*}
   Z(\ve_1,\ve_2;\q)
   = \sum_{n=0}^\infty \q^n Z_n(\ve_1,\ve_2)
   = \sum_{n=0}^\infty \q^n
   (\iota_{0*})^{-1} \pi_* [\operatorname{Hilb}^n\C^2],
\end{equation*}
where $\pi\colon\operatorname{Hilb}^n\C^2\to S^n\C^2$ is the
Hilbert-Chow morphism and $\iota_0$ is the inclusion of the unique
fixed point $n[0]$ in $S^n\C^2$. By \thmref{thm:M(r,n)weights} this is 
equal to \eqref{eq:formula1}.

Next we consider the Hilbert scheme
$\operatorname{Hilb}^n\widehat{\C}^2$ of $n$ points on the blowup
$\widehat{\C}^2$. The fixed points with respect to the
$\C^*\times\C^*$-action are parametrized by pairs of Young diagrams
$(Y^1,Y^2)$ by \propref{prop:fixedpoint2}.

Let $\mu(C)\in
H^2_{\C^*\times\C^*}(\operatorname{Hilb}^n\widehat{\C}^2)$ be the
class attached to the exceptional divisor $C$. (See the next section
for the definition.) We then define the partition function on the
blowup by
\begin{equation*}
   \widehat{Z}(\ve_1,\ve_2;t;\q)
   = \sum_{n=0}^\infty \q^n \sum_{d=0}^\infty \frac{t^d}{d!}
     \widehat{Z}_{n,d}(\ve_1,\ve_2)
   = \sum_{n=0}^\infty \q^n \sum_{d=0}^\infty \frac{t^d}{d!}
   (\iota_{0*})^{-1}\widehat\pi_*\left(\mu(C)^d\cap[
     \operatorname{Hilb}^n\widehat{\C}^2]\right),
\end{equation*}
where $\widehat\pi$ is the composite of the Hilbert-Chow morphism
$\operatorname{Hilb}^n\widehat{\C}^2\to S^n\widehat{\C}^2$ and 
the morphism $S^n\widehat{\C}^2\to S^n\C^2$.

By the \lemref{lem:mu} below, we have
\begin{equation*}
   \iota_{(Y^1,Y^2)}^* \mu(C) = |Y^1|\ve_1 + |Y^2|\ve_2.
\end{equation*}
Together with \thmref{thm:bM(r,k,n)weights} we have
\begin{equation*}
   \sum_{d=0}^\infty \frac{t^d}{d!} \sum_{(Y^1,Y^2)}
   \frac{(|Y^1|\ve_1 + |Y^2|\ve_2)^d q^{|Y^1| + |Y^2|}}
   {n_{Y^1}(\ve_1,\ve_2-\ve_1)\, n_{Y^2}(\ve_1-\ve_2,\ve_2)},
\end{equation*}
where $n_Y(\ve_1,\ve_2)$ is the denominator of
\eqref{eq:formula1}. This is equal to
\begin{equation*}
\begin{split}
   & \sum_{(Y^1,Y^2)} \frac{(\q e^{t\ve_1})^{|Y^1|}\, (\q e^{t\ve_2})^{|Y^2|}}
   {n_{Y^1}(\ve_1,\ve_2-\ve_1)\, n_{Y^2}(\ve_1-\ve_2,\ve_2)}
  = Z(\ve_1,\ve_2-\ve_1;\q e^{t\ve_1})Z(\ve_1-\ve_2,\ve_2;\q e^{t\ve_2})
\\
  =\; & \exp\left(\frac{\q e^{t\ve_1}}{\ve_1(\ve_2-\ve_1)} + 
  \frac{\q e^{t\ve_2}}{(\ve_1-\ve_2)\ve_2}\right).
\end{split}
\end{equation*}
We divide this by $Z(\ve_1,\ve_2;\q)$ and take the limit
$\ve_1,\ve_2\to 0$:
\begin{equation*}
   \lim_{\ve_1,\ve_2\to 0}
   \frac{\widehat{Z}(\ve_1,\ve_2;t;\q)}{Z(\ve_1,\ve_2;\q)}
   = \exp\left(-\frac{\q t^2}2\right).
\end{equation*}
This is a prototype of the blowup formula which will be discussed in
\secref{sec:blowupformula}. It should be noticed that this is equal to 
the generating function of 
\(
   \int_{\operatorname{Hilb}^n\widehat{X}} \mu(C)^{2n}
\)
for arbitrary smooth surface $X$. The minus sign comes from the
self-intersection number of $C$: $[C]^2 = -1$.
This can be shown roughly as follows: first show that $\mu(C)$ is a
pull-back of a class in $S^n\widehat{X}$ via the Hilbert-Chow morphism
$\operatorname{Hilb}^n\widehat{X}\to S^n\widehat{X}$. Then the
intersection numbers are those on $\widehat{X}^n$ divided by $n!$. The
class $\mu(C)$ corresponds to $\sum_i p_i^*[C]$, where $p_i\colon
\widehat{X}^n\to \widehat{X}$ is the $i$th projection.

\section{Instanton counting}\label{sec:counting}

We define {\it the partition function\/} as the following generating
function:
\begin{equation}\label{eq:defZ}
   Z(\ve_1,\ve_2,\vec{a};\q)
   = \sum_{n=0}^\infty \q^n Z_n(\ve_1,\ve_2,\vec{a})
   = \sum_{n=0}^\infty \q^n
   (\iota_{0*})^{-1} \pi_* [M(r,n)],
\end{equation}
where $[M(r,n)]$ denote the fundamental class of $H^{\hT}_*(M(r,n))$.
As we explained, this has an expression in terms of Hilbert
series~\eqref{eq:HilbPoly}. By (the equivariant homology analog of)
\propref{prop:Lef} together with \thmref{thm:M(r,n)weights}, we have
\begin{equation}\label{eq:Zsum}
    Z(\ve_1,\ve_2,\vec{a};\q)
   = \sum_{\vec{Y}}\frac{\q^{|\vec{Y}|}}{e(T_{\vec{Y}})}
   = \sum_{\vec{Y}} \frac{\q^{|\vec{Y}|}}
    {\displaystyle\prod_{\alpha,\beta}
    n^{\vec{Y}}_{\alpha,\beta}(\ve_1,\ve_2,\vec{a})},
\end{equation}
where
\begin{multline*}
n^{\vec{Y}}_{\alpha,\beta}(\ve_1,\ve_2,\vec{a})
   = \prod_{s \in Y_\alpha}
      \left( -l_{Y_\beta}(s)\ve_1 + (a_{Y_\alpha}(s)+1)\ve_2 + a_\beta 
   - a_\alpha\right)
\\
    \times\prod_{t\in Y_\beta} 
        \left((l_{Y_\alpha}(t)+1)\ve_1 -a_{Y_\beta}(t)\ve_2 + a_\beta 
   - a_\alpha\right)
.
\end{multline*}
This is nothing but Nekrasov's definition of the partition function
\cite[(1.6),(3.20)]{Nek}. We set
\begin{equation*}
   \Fin(\ve_1,\ve_2,\vec{a};\q)
   = \sum_{n=1}^\infty \q^n \Fin_n(\ve_1,\ve_2,\vec{a})
   = 
   \ve_1 \ve_2 \log Z(\ve_1,\ve_2,\vec{a};\q).
\end{equation*}

We give elementary properties of the partition function.
\begin{Lemma}\label{lem:elem}
\textup{(1)} $Z(\ve_2,\ve_1,\vec{a};\q) = Z(\ve_1,\ve_2,\vec{a};\q)$.

\textup{(2)} $Z(\ve_1,\ve_2, w\cdot\vec{a};\q) = Z(\ve_1,\ve_2,\vec{a};\q)$
where $w$ is an element of the symmetric group of $r$ letters.

\textup{(3)} $Z(-\ve_1,-\ve_2,-\vec{a};\q) = Z(\ve_1,\ve_2,\vec{a};\q)$.
\end{Lemma}

\begin{proof}
(1) The exchange of $\ve_1$ and $\ve_2$ is compensated with the
exchange of the Young diagram $Y_\alpha$ with its conjugate
$Y_\alpha'$.

(2) The exchange of $a_\alpha$ and $a_\beta$ is compensated with the
exchange of $Y_\alpha$ and $Y_\beta$.

(3) Clear from
\(
   n^{\vec{Y}}_{\alpha,\beta}(-\ve_1,-\ve_2,-\vec{a})
   = (-1)^{|Y_\alpha|+|Y_\beta|}
   n^{\vec{Y}}_{\alpha,\beta}(\ve_1,\ve_2,\vec{a}).
\)
\end{proof}

We now consider the moduli spaces on the blowup. By
\thmref{thm:bM(r,k,n)weights} the Euler class of the tangent space of
$\bM(r,k,n)$ at a fixed point $(\vec{k},\vec{Y}^1,\vec{Y}^2)$ is given
by
\begin{equation}\label{eq:Etanb}
   \prod_{\alpha,\beta} 
     l^{\vec{k}}_{\alpha,\beta}(\ve_1,\ve_2,\vec{a})
     \;
     n^{\vec{Y}^1}_{\alpha,\beta}(\ve_1,\ve_2-\ve_1,
       \vec{a} + \ve_1 \vec{k})
     \;
     n^{\vec{Y}^2}_{\alpha,\beta}(\ve_1-\ve_2,\ve_2,
        \vec{a} + \ve_2\vec{k}),
\end{equation}
where
\begin{equation*}
   l^{\vec{k}}_{\alpha,\beta}(\ve_1,\ve_2,\vec{a}) = 
   \begin{cases}
     {\displaystyle
     \prod_{\substack{i,j\ge 0\\i+j \le k_\alpha-k_\beta-1}}}
          (-i\ve_1 -j\ve_2 + a_\beta-a_\alpha)
       & \text{if $k_\alpha > k_\beta$}, \\
     {\displaystyle
     \prod_{\substack{i,j\ge 0\\i+j \le k_\beta-k_\alpha-2}}}
          \left((i+1)\ve_1 + (j+1)\ve_2 + a_\beta - a_\alpha\right)
       & \text{if $k_\alpha + 1 < k_\beta$}, \\
     1 & \text{otherwise}.
   \end{cases}
\end{equation*}
Note that $l^{\vec{k}}_{\alpha,\beta}(\ve_1,\ve_2,\vec{a})$ is
independent of $\vec{Y}^1$, $\vec{Y}^2$.

From now we use terminology for root
systems of Lie algebras as in \secref{sec:SWprep}, i.e.,
$\alpha_i\in\mathfrak h^*$, $\alpha_i^\vee\in\mathfrak h$,
$\vec{a} = \sum a^i\alpha_i^\vee$, etc.
Recall that $Q$ is the coroot lattice
\(
   \{ (k_1,\dots,k_r)\in \Z^r \mid \sum_\alpha k_\alpha = 0 \}.
\)
In order to treat the case $k = \sum_\alpha k_\alpha\neq 0$ (this
means that the gauge group is $\operatorname{\rm PU}(n)$ rather than
$\SU(n)$), we consider a normalization $\vec{l} = (k_1 - \frac{k}r,
\dots, k_r - \frac{k}r)$ as an element of the coweight lattice
\(
  P = \{ \vec{l} = (l_1,\dots,l_r)\in\Q^r \mid \sum_\alpha l_\alpha = 0,
  \exists k\in\Z\; \forall \alpha\; l_\alpha \equiv -\frac{k}r \mod \Z
  \}.
\)
There exists a homomorphism $P\to \Z/r\Z$ by taking the fractional
part of $l_\alpha$. It can be identified with the natural quotient
homomorphism $P\to P/Q$. We denote it by $\vec{l}\mapsto \{ \vec{l}\}$.
Hereafter we identify $\vec{l}$ with $\vec{k}$ and denote both by
$\vec{k}$. We write $\vec{k} = \sum_i k^i \alpha_i^\vee$ in either
case $k=0$, $\neq 0$. But $k^i$ may be rational in the latter case.
Let $(\ ,\ )$ be the standard inner product on $\mathfrak h$. The
Killing form $B_{\SU(r)}$ of $\SU(r)$ satisfies $B_{\SU(r)} = 2r(\ ,\ 
)$.
The following formulas are useful later:
\begin{equation}
\begin{gathered}
   \frac1{2r}   
   \sum_{\alpha, \beta} (k_\alpha - k_\beta)(a_\alpha - a_\beta)
   =
   (\vec{k},\vec{a})
   = \sum_{ij} C_{ij} a^i k^j,
\\
    \frac1{2r}\sum_{\alpha, \beta} (k_\alpha - k_\beta)^2
    =
    (\vec{k}, \vec{k})
    = \sum_{i,j} C_{ij} k^i k^j,
\\
  \sum_{\alpha < \beta} \frac{k_\alpha-k_\beta}2
  =
  \langle\vec{k},\rho\rangle
  = \sum_i k^i.
\end{gathered}
\label{eq:useful}\end{equation}
Here $C_{ij}$ is the Cartan matrix, and
$\rho$ is the half of the sum of positive roots, as usual.

For a root $\alpha\in\Delta$, we define
\begin{equation}\label{eq:l}
   l^{\vec{k}}_{\alpha}(\ve_1,\ve_2,\vec{a}) = 
   \begin{cases}
     {\displaystyle
     \prod_{\substack{i,j\ge 0\\i+j \le -\langle\vec{k}, \alpha\rangle-1}}}
          (-i\ve_1 -j\ve_2 + \langle \vec{a}, \alpha\rangle)
       & \text{if $\langle \vec{k}, \alpha\rangle < 0$}, \\
     {\displaystyle
     \prod_{\substack{i,j\ge 0\\i+j\le \langle \vec{k}, \alpha\rangle-2}}}
          \left((i+1)\ve_1 + (j+1)\ve_2 + \langle \vec{a},\alpha\rangle\right)
       & \text{if $\langle\vec{k}, \alpha\rangle > 1$}, \\
     1 & \text{otherwise},
   \end{cases}
\end{equation}
where $l^{\vec{k}}_{\alpha,\beta}$ in the previous notation
corresponds to $l^{\vec{k}}_{e_{\beta,\alpha}}$.

The following will be useful later:
\begin{Lemma}\label{lem:lsym}
\textup{(1)}
\(
   l^{\vec{k}}_{\alpha}(\ve_1,\ve_2,\vec{a}) = 
   l^{\vec{k}}_{\alpha}(\ve_2,\ve_1,\vec{a}).
\)

\textup{(2)}
\(
   l^{\vec{k}}_{\alpha}(\ve_1,\ve_2,\vec{a}) = 
   (-1)^{\langle\vec{k},\alpha\rangle(\langle\vec{k},\alpha\rangle-1)/2}\,
   l^{-\vec{k}}_{-\alpha}(-\ve_1,-\ve_2,\vec{a}).
\)

\textup{(3)}
$l^{\vec{k}}_{\alpha}(\ve_1,\ve_2,\vec{a})$ is regular at
$(\ve_1,\ve_2) = 0$ and
\begin{equation*}
   l^{\vec{k}}_{\alpha}(0,0,\vec{a}) = 
     \langle \vec{a},\alpha\rangle^{
     \langle\vec{k},\alpha\rangle(\langle\vec{k},\alpha\rangle-1)/2}
\end{equation*}
\end{Lemma}

Let $\mathcal E$ be a universal sheaf on $\bp\times\bM(r,k,n)$. We
define an equivariant cohomology class $\mu(C)\in
H^2_{\hT}(\bM(r,k,n))$ by
\begin{equation*}
   (c_2(\mathcal E)-\frac{r-1}{2r}c_1(\mathcal E)^2)/[C],
\end{equation*}
where $/$ denotes the slant product
\(
  /\colon H^d_{\hT}(\bp\times\bM(r,k,n))\otimes H_i^{\hT}(\bp)
  \to H^{d-i}_{\hT}(\bM(r,k,n)).
\)
Note that we have
\(
   c_2(\mathcal E)-\frac{r-1}{2r}c_1(\mathcal E)^2
   = \frac1{2r}c_2(\End\mathcal E)
\)
on the open locus $\bM^{\operatorname{reg}}(r,k,n)$.

Let $\iota_{(\vec{k},\vec{Y}^1,\vec{Y}^2)}$ be the inclusion of the
fixed point $(\vec{k},\vec{Y}^1,\vec{Y}^2)$ into $\bM(r,k,n)$.
\begin{Lemma}\label{lem:mu}
\begin{equation*}
   \iota^*_{(\vec{k},\vec{Y}^1,\vec{Y}^2)}\mu(C)
   = |\vec{Y}^1|\,\ve_1 + |\vec{Y}^2|\,\ve_2
   + (\vec{k},\vec{a})
   + \frac{(\vec{k},\vec{k})}2(\ve_1+\ve_2)
.
\end{equation*}
\end{Lemma}

\begin{proof}
Let $E$ be a sheaf corresponding to the fixed point
$(\vec{k},\vec{Y}^1,\vec{Y}^2)$. We have
\begin{equation*}
   c_2(E)-\frac{r-1}{2r}c_1(E)^2 
   =  |\vec{Y}^1|\, [p_1] + |\vec{Y}^2|\, [p_2]
      + c_2(E^{\vee\vee})-\frac{r-1}{2r}c_1(E^{\vee\vee})^2.
\end{equation*}
The double dual $E^{\vee\vee}$ is a direct sum $\bigoplus_\alpha {\cal
O}_X(k_\alpha C)e_\alpha$. Therefore
\begin{equation*}
   c_2(E^{\vee\vee})-\frac{r-1}{2r}c_1(E^{\vee\vee})^2
   = \frac1{2r}c_2(\End E^{\vee\vee})
   = -\frac{1}{2r}\sum_{\alpha<\beta}
   \left\{(k_\alpha[C]+a_\alpha)-(k_\beta[C]+a_\beta)\right\}^2.
\end{equation*}
Substituting
\[
   \int_{\bp} [p_1][C] = \ve_1,
\quad
   \int_{\bp} [p_2][C] = \ve_2,
\quad
   \int_{\bp} [C]=0, 
\quad
   \int_{\bp} [C]^2=-1,
\quad
   \int_{\bp} [C]^3=-(\ve_1+\ve_2),
\]
into this, we get
\begin{equation*}
   \iota^*_{(\vec{k},\vec{Y}^1,\vec{Y}^2)}\mu(C)
   = |\vec{Y}^1|\,\ve_1 + |\vec{Y}^2|\,\ve_2
   + \frac{1}{2r}\sum_{\alpha<\beta}
   (2(k_\alpha-k_\beta)(a_\alpha-a_\beta)+(k_\alpha-k_\beta)^2(\ve_1+\ve_2)).
\end{equation*}
This is the desired formula thanks to \eqref{eq:useful}.
\end{proof}

We now define the partition function on the blowup:
\begin{equation*}
   \widehat{Z}^{k}(\ve_1,\ve_2,\vec{a};t;\q)
   = \sum_{n} \q^n \sum_{d=0}^\infty \frac{t^d}{d!}
     \widehat{Z}^{k}_{n,d}(\ve_1,\ve_2,\vec{a})
   = \sum_{n} \q^n \sum_{d=0}^\infty \frac{t^d}{d!}
   (\iota_{0*})^{-1}\widehat\pi_*\left(\mu(C)^d\cap[{\bM(r,k,n)}]\right),
\end{equation*}
where $n$ runs over $\Z_{\ge 0} - \frac1{2r}k(r-k)$.
By (\ref{eq:Etanb}, \ref{lem:mu}) this can be represented in terms of
Nekrasov's partition function:
{\allowdisplaybreaks
\begin{equation}\label{eq:Zhat}
\begin{split}
  & \widehat{Z}^{k}_{n,d}(\ve_1,\ve_2,\vec{a})
  =
  \sum_{\frac12(\vec{k},\vec{k})+l+m = n}
  \begin{aligned}[t]
   & \left(
     l\ve_1 + m\ve_2
   + (\vec{k}, \vec{a})
   + \frac{(\vec{k},\vec{k})}2(\ve_1+\ve_2)
   \right)^d \frac1
   {\prod_{\alpha\in\Delta} l^{\vec{k}}_{\alpha}(\ve_1,\ve_2,\vec{a})}
   \\
   &\qquad
   \times
   Z_l(\ve_1,\ve_2-\ve_1,\vec{a}+\ve_1\vec{k})
   Z_m(\ve_1-\ve_2,\ve_2,\vec{a}+\ve_2\vec{k})
.
  \end{aligned}
\end{split}
\end{equation}
The generating function is} 
{\allowdisplaybreaks
\begin{equation}\label{eq:blowup}
\begin{split}
   &
   \widehat{Z}^{k}(\ve_1,\ve_2,\vec{a};t;\q)
   =
  \begin{aligned}[t]
  & \sum_{\{\vec{k}\} = -\frac{k}r}
   \exp\left[t\left(
   (\vec{k}, \vec{a})
   + \frac{(\vec{k},\vec{k})}2(\ve_1+\ve_2)
   \right)\right]
   \frac{\q^{\frac12(\vec{k},\vec{k})}}
   {\prod_{\alpha\in\Delta} l^{\vec{k}}_{\alpha}(\ve_1,\ve_2,\vec{a})}
\\     
   &\qquad\times
   Z(\ve_1,\ve_2-\ve_1,\vec{a}+\ve_1\vec{k};\q e^{t\ve_1})\,
   Z(\ve_1-\ve_2,\ve_2,\vec{a}+\ve_2\vec{k};\q e^{t\ve_2})
.
  \end{aligned}
\end{split}
\end{equation}
}

We now only consider the case $k = 0$ for a while. We omit the
superscript in this case.

\begin{Proposition}\label{prop:lowblowup}
\textup{(1)} 
\(
   \widehat\pi_*[{\bM(r,0,n)}] = [M_0(r,n)].
\)

\textup{(2)} 
\(
   \widehat\pi_*\left(\mu(C)^d\cap[{\bM(r,0,n)}]\right) = 0
\)
for $1 \le d \le 2r-1$.
\end{Proposition}

\begin{proof}
Both results are well-known in Donaldson theory (see e.g.,
\cite[3.8.1]{FM}). We give a proof for the completeness.

(1) By the dimension reason, the inclusion $i$ of
$M_0^{\operatorname{reg}}(r,n)$ in $M_0(r,n)$ induces an isomorphism
in degree $4nr$:
\begin{equation*}
   H^{\hT}_{4nr}(M_0(r,n)) \xrightarrow[\cong]{i^*}
   H^{\hT}_{4nr}(M_0^{\operatorname{reg}}(r,n)).
\end{equation*}
Therefore it is enough to show that 
\(
   i^*\widehat\pi_*([\bM(r,0,n)]) = [M_0^{\operatorname{reg}}(r,n)].
\)
But this is clear since $\widehat\pi$ becomes an isomorphism over the set
$M_0^{\operatorname{reg}}(r,n)$.

(2) First note that $\mu(C)$ is equal to $c_1(\mathcal L)$, where
$\mathcal L$ is the determinant line bundle over $\bM(r,0,n)$ where
the fiber over $(E,\Phi)$ is
\begin{equation*}
   \left(\Lambda^{\max} H^1(\bp, E(-\linf))\right)^*
   \otimes \Lambda^{\max} H^1(\bp, E(C-\linf)).
\end{equation*}
This line bundle has a natural section $s$ whose zero set is a
representative of $\mu(C) = c_1(\mathcal L)$ and consists of bundles
that restrict to $C$ in a non-trivial way. (See \cite[4.6]{Bryan:1997}.)

Consider
\[
   \overline{\{ 0\} \times M_0^{\operatorname{reg}}(r,n-1)
   }.
\]
This has complex codimension $2r$. Therefore if $i\colon U\to
M_0(r,n)$ denote the inclusion of the complement, the pullback
homomorphism $i^*$ is an isomorphism in degree $\ge 4nr - 4r + 2$.
Therefore we can restrict $\widehat\pi$ to $U$ as in (1). Now the
vanishing is clear since the section $s$ of $\mathcal L$ does not
vanish there as explained above.
\end{proof}

If we apply this to \eqref{eq:Zhat} with $d=1,2$, we get
\begin{equation*}
\begin{split}
  n\ve_1 Z_n(\ve_1,\ve_2-\ve_1,\vec{a}) +
  n\ve_2 Z_n(\ve_1-\ve_2,\ve_2,\vec{a}) &= A,
\\
  n^2\ve_1^2 Z_n(\ve_1,\ve_2-\ve_1,\vec{a}) +
  n^2\ve_2^2 Z_n(\ve_1-\ve_2,\ve_2,\vec{a}) &= B,
\end{split}
\end{equation*}
where $A$ and $B$ are given by lower terms $l,m < n$ (but $\vec{a}$
may be shifted by $\vec{k}$). Therefore we can determine
$Z_n(\ve_1,\ve_2-\ve_1,\vec{a})$, $Z_n(\ve_1-\ve_2,\ve_2,\vec{a})$
recursively. Changing $\ve_2$ by $\ve_1+\ve_2$, we get
$Z_n(\ve_1,\ve_2, \vec{a})$.
%

In order to express these assertions by differential equations, we
introduce the following generalization of the Hirota differential:
\begin{equation*}
\begin{split}
   \left(D^{(\ve_1,\ve_2)}_x\right)^m (f\cdot g)
   &= \left.(\frac{d}{dy})^m f(x + \ve_1 y) g(x + \ve_2 y)\right|_{y = 0}
\\
   &= \sum_{k=0}^m \ve_1^k \ve_2^{m-k} {m\choose k}
   \frac{d^k f}{dx^k}\frac{d^{m-k}g}{dx^{m-k}}. 
\end{split}
\end{equation*}
$(D^{(1,-1)}_x)^m$ is the ordinary Hirota differential.
We have
\begin{equation*}
\begin{split}
   & Z(\ve_1,\ve_2-\ve_1,\vec{a}+\ve_1\vec{k};\q e^{t\ve_1})
   Z(\ve_1-\ve_2,\ve_2,\vec{a}+\ve_2\vec{k};\q e^{t\ve_2})
\\
   =\; &
   \exp(tD^{(\ve_1,\ve_2)}_{\log \q})
   \left(Z(\ve_1,\ve_2-\ve_1,\vec{a}+\ve_1\vec{k};\q)\cdot
   Z(\ve_1-\ve_2,\ve_2,\vec{a}+\ve_2\vec{k};\q)\right).
\end{split}
\end{equation*}

\begin{Corollary}
The followings hold:
{\allowdisplaybreaks
\begin{gather}
  Z(\ve_1,\ve_2,\vec{a};\q) = 
  \sum_{\vec{k}}
   \frac{\q^{\frac12(\vec{k},\vec{k})}}
   {\prod_{\alpha\in\Delta} l^{\vec{k}}_{\alpha}(\ve_1,\ve_2,\vec{a})}
   Z(\ve_1,\ve_2-\ve_1,\vec{a}+\ve_1\vec{k};\q)
   Z(\ve_1-\ve_2,\ve_2,\vec{a}+\ve_2\vec{k};\q)
   \label{eq:Zind}
\\
  0 =
  \sum_{\vec{k}}
   \frac{\q^{\frac12(\vec{k},\vec{k})}}
   {\prod_{\alpha\in\Delta} l^{\vec{k}}_{\alpha}(\ve_1,\ve_2,\vec{a})}
   \begin{aligned}[t]
   & \left(D^{(\ve_1,\ve_2)}_{\log \q}
   + (\vec{k}, \vec{a})
   + \frac{(\vec{k},\vec{k})}2(\ve_1+\ve_2)
   \right)^d
 \\
   &\qquad \left(Z(\ve_1,\ve_2-\ve_1,\vec{a}+\ve_1\vec{k};\q)\cdot
   Z(\ve_1-\ve_2,\ve_2,\vec{a}+\ve_2\vec{k};\q)\right)
   \end{aligned}
  \label{eq:ind}
\end{gather}
for $1\le d\le 2r-1$.}
\end{Corollary}

The second equation \eqref{eq:ind} will play a fundamental role in our
study of the partition function $Z$. We call it the {\it blowup
  equation}.

\begin{proof}
\eqref{eq:Zind} follows from \propref{prop:lowblowup}(1) by setting
$t=0$ in \eqref{eq:Zhat}. \propref{prop:lowblowup}(2) means that
$(\frac{d}{dt})^d\widehat{Z}(\ve_1,\ve_2,\vec{a};t;q)|_{t=0} = 0$ with
$1\le d \le 2r-1$. We get the above if we differentiate the right hand
side of \eqref{eq:Zhat}.
\end{proof}

For a later purpose, we divide the blowup equations for $d=1,2$ by
$Z(\ve_1,\ve_2-\ve_1,\vec{a};\q)Z(\ve_1-\ve_2,\ve_2,\vec{a};\q)$
and write down explicitly as
{\allowdisplaybreaks
\begin{gather}
  0 =
  \begin{aligned}[t]
  & \sum_{\vec{k}}
   \frac{\q^{\frac12(\vec{k},\vec{k})}}
   {\prod_{\alpha\in\Delta} l^{\vec{k}}_{\alpha}(\ve_1,\ve_2,\vec{a})}
  \Biggl[
     \frac1{\ve_2 - \ve_1}\left(
     \q\frac{\partial}{\partial \q}\Fin_a(\vec{a}+\ve_1\vec{k})
     -
     \q\frac{\partial}{\partial \q}\Fin_b(\vec{a}+\ve_2\vec{k})
     \right)
\\
   & \qquad\qquad\qquad\qquad\qquad\qquad
     + (\vec{k}, \vec{a})
     + \frac{(\vec{k},\vec{k})}2(\ve_1+\ve_2)
     \Biggr]
\\     
   &
   \times\exp\left[
   \frac1{\ve_2-\ve_1}
   \left(
   \frac{ \Fin_a(\vec{a}+\ve_1\vec{k}) - \Fin_a(\vec{a}) }{\ve_1}
   -
   \frac{ \Fin_b(\vec{a}+\ve_2\vec{k}) - \Fin_b(\vec{a}) }{\ve_2}
    \right)\right],
  \end{aligned}
  \label{eq:ind1}
\\
  0 =
  \begin{aligned}[t]
  & \sum_{\vec{k}}
   \frac{\q^{\frac12(\vec{k},\vec{k})}}
   {\prod_{\alpha\in\Delta} l^{\vec{k}}_{\alpha}(\ve_1,\ve_2,\vec{a})}
   \Biggl[\biggl\{
     (\vec{k}, \vec{a})
     + \frac{(\vec{k},\vec{k})}2(\ve_1+\ve_2)
\\
   &\qquad\qquad\qquad\qquad\qquad
     + \frac1{\ve_2-\ve_1}\left(
       \q\frac{\partial}{\partial \q}\Fin_a(\vec{a}+\ve_1\vec{k})
       - \q\frac{\partial}{\partial \q}\Fin_b(\vec{a}+\ve_2\vec{k})\right)
       \biggr\}^2
\\
   &\qquad\qquad\qquad
    + \frac1{\ve_2-\ve_1}\left(
     {\ve_1}(\q\frac{\partial}{\partial \q})^2 \Fin_a(\vec{a}+\ve_1\vec{k})
     - \ve_2(\q\frac{\partial}{\partial \q})^2 \Fin_b(\vec{a}+\ve_2\vec{k})
   \right)
   \Biggr]
\\
   &
   \times\exp\left[
   \frac1{\ve_2-\ve_1}
   \left(
   \frac{ \Fin_a(\vec{a}+\ve_1\vec{k}) - \Fin_a(\vec{a}) }{\ve_1}
   -
   \frac{ \Fin_b(\vec{a}+\ve_2\vec{k}) - \Fin_b(\vec{a}) }{\ve_2}
    \right)\right],
   \end{aligned}\label{eq:ind2}
\end{gather}
where}
\begin{gather*}
  \exp\frac{\Fin_a(\vec{a})}{\ve_1(\ve_2-\ve_1)}
  = Z(\ve_1,\ve_2-\ve_1,\vec{a};\q),
  \qquad
  \exp\frac{\Fin_b(\vec{a})}{(\ve_1-\ve_2)\ve_2}
  = Z(\ve_1-\ve_2,\ve_2,\vec{a};\q).
\end{gather*}
The functions $\Fin_a$, $\Fin_b$ depends also on $\ve_1$, $\ve_2$, but
we omit them from the notation for brevity.

\section{Behavior at $\ve_1, \ve_2 = 0$}\label{sec:limit}

We prove Nekrasov's conjecture in this section.

\begin{Lemma}\label{lem:symZ}
$Z(\ve_1,-2\ve_1,\vec{a};\q) = Z(2\ve_1, -\ve_1,\vec{a};\q)$.
\end{Lemma}

\begin{proof}
We set $\ve_2 = -\ve_1$ in \eqref{eq:Zhat} with $d=1$. Then we have
\begin{equation*}
\begin{split}
  & n\ve_1 \left( Z_n(\ve_1,-2\ve_1,\vec{a}) -
  Z_n(2\ve_1,-\ve_1,\vec{a})\right)
\\
  = \; &
  \sum_{\substack{\frac12(\vec{k},\vec{k})+l+m = n\\
      l\neq n, m\neq n}}
   \left\{ (l - m)\ve_1 + (\vec{k},\vec{a}) \right\}
   \;
   \frac{
   Z_l(\ve_1,-2\ve_1,\vec{a}+\ve_1\vec{k})
   Z_m(2\ve_1,-\ve_1,\vec{a}-\ve_1\vec{k})}
   {\prod_{\alpha\in\Delta} l^{\vec{k}}_{\alpha}(\ve_1,-\ve_1,\vec{a})}
.
\end{split}
\end{equation*}
We show
\(
  Z_n(\ve_1,-2\ve_1,\vec{a}) = Z_n(2\ve_1, -\ve_1,\vec{a})
\)
by the induction on $n$. The assertion is trivial for $n=1$. Suppose
that it is true for $l,m< n$. Then the right hand side of the above
equation vanishes, as terms with $(\vec{k},l,m)$ and $(-\vec{k}, m,
l)$ cancel and the term $(0,l,l)$ is $0$.
Here we have used
\(
   l^{\vec{k}}_{\alpha}(\ve_1,-\ve_1,\vec{a})
   = (-1)^{\langle\vec{k},\alpha\rangle(\langle\vec{k},\alpha\rangle-1)/2}
   l^{-\vec{k}}_{-\alpha}(\ve_1,-\ve_1,\vec{a})
\)
which follows from \lemref{lem:lsym}, and that
\[
    \sum_{\alpha\in\Delta}
    \langle \vec{k},\alpha\rangle(\langle \vec{k},\alpha\rangle - 1)/2
    = r (\vec{k},\vec{k})
\]
is an even number.
\end{proof}

The following follows from this lemma and its proof:
\begin{Corollary}
$\widehat{Z}_{n,d}(\ve_1,-\ve_1,\vec{a})$ vanishes for odd $d$.
\end{Corollary}

This is compatible with what is known for the usual blowup formula for 
Donaldson invariants (cf.\ \cite{Fintushel-Stern:1996}.)

The following is the first part of Nekrasov's conjecture:
\begin{Proposition}
$\Fin(\ve_1,\ve_2,\vec{a};\q)$ is regular at $\ve_1 = \ve_2 = 0$.
\end{Proposition}

\begin{proof}
The point of the proof is a recursive structure of the blowup equation
(\ref{eq:ind1}, \ref{eq:ind2}). Let us separate terms with $\vec{k} =
0$:
\begin{equation*}
\begin{split}
    & \frac1{\ve_2 - \ve_1}\left(
     \q\frac{\partial}{\partial \q}\Fin_a(\vec{a})
     -
     \q\frac{\partial}{\partial \q}\Fin_b(\vec{a})
     \right)
   = A,
\\
    &
    \begin{aligned}[b]
    & \frac1{(\ve_2-\ve_1)^2}\left(
       \q\frac{\partial}{\partial \q}\Fin_a(\vec{a})
       - \q\frac{\partial}{\partial \q}\Fin_b(\vec{a})\right)^2
     \\
     &\qquad
     + \frac1{\ve_2-\ve_1}\left(
     {\ve_1}(\q\frac{\partial}{\partial \q})^2 \Fin_a(\vec{a})
     - \ve_2(\q\frac{\partial}{\partial \q})^2 \Fin_b(\vec{a})\right)
    \end{aligned}
    = B,
\end{split}
\end{equation*}
where $A$ and $B$ are terms with $\vec{k}\neq 0$ and hence divisible
by $\q$. We further replace the first term in the second equation by
$A^2$. If we express $\Fin_a(\vec{a})$, $\Fin_b(\vec{a})$ by formal
power series in $\q$, then the above equations determine the
coefficients recursively. We want to show that
$\Fin_n(\ve_1,\ve_2,\vec{a})$ is regular at $\ve_1 = \ve_2 = 0$ by the
induction using this recursive system.
This is equivalent to showing that $A$ and $B$ are regular under the
assumption that $\Fin_a(\vec{a})$, $\Fin_b(\vec{a})$ are regular. This 
follows from the following lemma.
\end{proof}

\begin{Lemma}\label{lem:limit}
Suppose that $\Fin(\ve_1,\ve_2,\vec{a};\q)$ is regular at
$(\ve_1,\ve_2) = 0$. Then the following are also regular and their
values are given by
\begin{gather*}
   \left.
    \frac1{\ve_2 - \ve_1}\left(
     \q\frac{\partial}{\partial \q}\Fin_a(\vec{a}+\ve_1\vec{k})
     -
     \q\frac{\partial}{\partial \q}\Fin_b(\vec{a}+\ve_2\vec{k})
     \right)\right|_{(\ve_1,\ve_2) = 0}
   =
   - \sum_i k^i\,
    \q\frac{\partial^2\Fin}{\partial \q\partial a^i}(0,0,\vec{a};\q),
\\
  \left.
  \frac1{\ve_2-\ve_1}\left(
     {\ve_1}(\q\frac{\partial}{\partial \q})^2 \Fin_a(\vec{a})
     - \ve_2(\q\frac{\partial}{\partial \q})^2 \Fin_b(\vec{a})\right)
   \right|_{(\ve_1,\ve_2) = 0}
  =
     - (\q\frac{\partial}{\partial\q})^2 \Fin(0,0,\vec{a};\q),
\\
\begin{aligned}[t]
   & \left.
   \frac1{\ve_2-\ve_1}
   \left(
   \frac{ \Fin_a(\vec{a}+\ve_1\vec{k}) - \Fin_a(\vec{a}) }{\ve_1}
   -
   \frac{ \Fin_b(\vec{a}+\ve_2\vec{k}) - \Fin_b(\vec{a}) }{\ve_2}
    \right)
  \right|_{(\ve_1,\ve_2) = 0}
\\
 & \qquad\qquad\qquad
  = - \frac12 \sum_{i,j}
  \frac{\partial^2\Fin}{\partial a^i\partial a^j}(0,0,\vec{a};\q) k^i k^j.
\end{aligned}
\end{gather*}
\end{Lemma}

\begin{proof}
The regularity is a consequence of the symmetry
\( 
   \Fin_b(\vec{a}) =
   \left.\Fin_a(\vec{a})\right|_{\ve_1\leftrightarrow\ve_2}.
\)
In order to show the above equalities, we just need to note
\[
   \frac{\partial\Fin}{\partial \ve_1}(0,0,\vec{a};\q) = 
   \frac{\partial\Fin}{\partial \ve_2}(0,0,\vec{a};\q) = 0.
\]
The first equality is the consequence of \lemref{lem:elem}(1),
and the second equality follows from \lemref{lem:symZ}.
\end{proof}

We now take the limit $\ve_1,\ve_2 \to 0$ of \eqref{eq:ind2}. (The
limit of \eqref{eq:ind1} becomes the trivial identity $0 = 0$.)
We set
\(
   \Finz(\vec{a};q) = \Fin(0,0,\vec{a};\q).
\)
By \lemref{lem:lsym}(3), we have
\begin{equation*}
   \prod_{\alpha\in\Delta}l^{\vec{k}}_{\alpha}(0,0,\vec{a})
   = \prod_{\alpha\in\Delta^+} 
    (-1)^{\langle \vec{k},\alpha\rangle(\langle\vec{k},\alpha\rangle+1)/2}\,
    \langle\vec{a},\alpha\rangle^{\langle\vec{k},\alpha\rangle^2}
   =
   (-1)^{\langle\vec{k},\rho\rangle}\,
   \prod_{\alpha\in\Delta^+} 
    \left\{\sqrt{-1}\langle \vec{a},\alpha\rangle\right\}^{
    \langle\vec{k},\alpha\rangle^2}
   .
\end{equation*}
Therefore we get
\begin{equation}\label{eq:Wh}
  0 =
  \begin{aligned}[t]
  & \sum_{\vec{k}}
   \frac
   {(-1)^{\langle\vec{k},\rho\rangle}\q^{\frac12(\vec{k},\vec{k})}}
   {
   \prod_{\alpha\in\Delta^+} 
    \left\{\sqrt{-1}\langle \vec{a},\alpha\rangle\right\}^{
    \langle\vec{k},\alpha\rangle^2}
    }
\\
   &\qquad\qquad
   \times \left[\left\{
     \sum_i k^i \left(
       \sum_j C_{ij} a^j
       -
     \q\frac{\partial^2\Finz}{\partial \q\partial a^i}(\vec{a};\q)
     \right)
       \right\}^2
   - (\q\frac{\partial}{\partial\q})^2 \Finz(\vec{a};\q)
   \right]
\\
   &\qquad\qquad\qquad\qquad
   \times
   \exp\left(
     - \frac12 \sum_{i,j}
     \frac{\partial^2\Finz}{\partial a^i\partial a^j}(\vec{a};\q) k^i k^j
   \right).
   \end{aligned}
\end{equation}
In order to compare this with the formula in literature, we introduce
the following functions:
\begin{gather}
   \tau_{ij} = 
   \frac{\sqrt{-1}}{\pi} \sum_{\alpha\in\Delta_+}
   \langle \alpha_i^\vee, \alpha\rangle 
   \langle \alpha_j^\vee, \alpha\rangle 
     \log\left(\frac{\sqrt{-1}\langle \vec{a},\alpha\rangle}
       {\q^{\frac1{2r}}}\right)
    -
    \frac1{2\pi\sqrt{-1}}
    \frac{\partial^2\Finz}{\partial a^i\partial a^j}(\vec{a};\q),
    \label{eq:tau'}
\\
   u_2 = \frac12 (\vec{a},\vec{a}) -
   \q\frac{\partial\Finz}{\partial \q}(\vec{a};\q).
   \label{eq:RG'}
\end{gather}
Now \eqref{eq:Wh} can be written as
\begin{equation}\label{eq:Wh2}
   (\q\frac{\partial}{\partial\q})^2 \Finz(\vec{a};\q)
   = \sum_{i,j} \frac{\partial u_2}{\partial a^i}
   \frac{\partial u_2}{\partial a^j}
   \frac1{\pi\sqrt{-1}}
   \frac{\partial}{\partial\tau_{ij}} \log\Theta_E(0 | \tau),
\end{equation}
where we have used \eqref{eq:useful} several times and $\Theta_E$ is
as in \eqref{eq:Theta}.
This, combined with \eqref{eq:RG'}, is exactly the contact term
equation \eqref{eq:contactterm} if we replace $\q^{\frac1{2r}}$ by
$\Lambda$. Note also that \eqref{eq:tau'} coincides with \eqref{eq:tau}.
And \eqref{eq:RG'} is nothing but \eqref{eq:RG}.

The equation \eqref{eq:Wh2} has the same structure as the blowup
equation \eqref{eq:ind}. When we expand $\Finz$ as a formal power
series in $\q$, coefficients are determined inductively. In
particular, the solution to the above equation is {\it unique}. This
observation was due to \cite{EMM}. (See also \cite{Matone} for an
earlier result for $\SU(2)$.) Since the Seiberg-Witten prepotential
satisfies \eqref{eq:Wh2}, we conclude that $\Finz$ coincides with its
instanton part. This is our confirmation of Nekrasov's conjecture.

\section{Blowup formula}\label{sec:blowupformula}

We divide \eqref{eq:blowup} by
$Z(\ve_1,\ve_2-\ve_1,\vec{a};\q)Z(\ve_1-\ve_2,\ve_2,\vec{a};\q)$ and
take the limit $\ve_1, \ve_2 \to 0$. We need the following
generalization of the third equation in \lemref{lem:limit}:
\begin{equation*}
\begin{split}
   &
\begin{aligned}[t]
   &
   \frac1{\ve_2-\ve_1}
   \biggl(
   \frac{
   \Fin(\ve_1,\ve_2-\ve_1,\vec{a}+\ve_1\vec{k};\q e^{t\ve_1})
    - \Fin(\ve_1,\ve_2-\ve_1,\vec{a};\q)
    }{\ve_1}
\\
   &\qquad\qquad\qquad\qquad
   -
   \frac{
   \Fin(\ve_1-\ve_2,\ve_2,\vec{a}+\ve_2\vec{k};\q e^{t\ve_2})
    - \Fin(\ve_1-\ve_2,\ve_2,\vec{a};\q)
    }{\ve_2}
    \biggr)
  \Biggr|_{(\ve_1,\ve_2) = 0}
  \end{aligned}
\\
 = \; &
  - \frac12 \left(\q \frac{\partial}{\partial\q}\right)^2
           \Finz(\vec{a};\q)\, t^2
  - \frac12 \sum_{i,j}
  \frac{\partial^2\Finz}{\partial a^i\partial a^j}(\vec{a};\q) k^i k^j
  - \sum_i \q\frac{\partial^2\Finz}{\partial\q\partial a^i}(\vec{a};\q)
  \, t k^i
  .
\end{split}
\end{equation*}
Here we have used
\(
  \frac{\partial}{\partial\log\q} = \q \frac{\partial}{\partial\q}.
\)
Therefore we get
\begin{Theorem}
$\widehat{Z}^{k}(\ve_1,\ve_2,\vec{a};t;\q)/Z(\ve_1,\ve_2,\vec{a};q)$ is
regular at $(\ve_1,\ve_2) = 0$. Its value is
\begin{multline*}
    \exp\left( - \frac12 \left(\q \frac{\partial}{\partial\q}\right)^2
           \Finz(\vec{a};\q)\, t^2 \right)
\\
    \times \sum_{\vec{k}\in P : \{ \vec{k}\} = -\frac{k}r} \Biggl\{
    \begin{aligned}[t]
    & \frac
   {(-1)^{\langle\vec{k},\rho\rangle}\q^{\frac12(\vec{k},\vec{k})}}
   {
   \prod_{\alpha\in\Delta^+} 
    \left\{\sqrt{-1}\langle \vec{a},\alpha\rangle\right\}^{(\vec{k},\alpha)^2}
    }
\\
   &\quad
    \exp\left( - \frac12 \sum_{i,j}
  \frac{\partial^2\Finz}{\partial a^i\partial a^j}(\vec{a};\q) k^i k^j
  + \sum_i \left(
    \sum_j C_{ij} a^j
    -  \q\frac{\partial^2\Finz}{\partial\q\partial a^i}(\vec{a};\q) \right)
  \, t k^i \right)\Biggr\}
    \end{aligned}
\\
  \times \Biggl[
    \sum_{\vec{k}\in Q}
    \begin{aligned}[t]
    & \frac
   {(-1)^{\langle\vec{k},\rho\rangle}\q^{\frac12(\vec{k},\vec{k})}}
   {
   \prod_{\alpha\in\Delta^+} 
    \left\{\sqrt{-1}\langle \vec{a},\alpha\rangle\right\}^{(\vec{k},\alpha)^2}
    }
    \exp\left( - \frac12 \sum_{i,j}
  \frac{\partial^2\Finz}{\partial a^i\partial a^j}(\vec{a};\q)
  k^i k^j \right)
    \Biggr]^{-1}.
    \end{aligned}
\end{multline*}
\end{Theorem}

If we use the theta function in \eqref{eq:Theta}, this can be
written simply as
\begin{equation*}
   \exp\left( - \frac12 \left(\q \frac{\partial}{\partial\q}\right)^2
      \Finz(\vec{a};\q)\, t^2 \right)
   \frac{\Theta_k(\vec{\xi}|\tau)}{\Theta_E(0|\tau)},
  \quad\text{where }
  \xi^i = \frac{t}{2\pi\sqrt{-1}}
   \frac{\partial u_2}{\partial a^i}.
\end{equation*}
Here $\Theta_k$ is defined as in \eqref{eq:Theta} where the summation
is over $\vec{k}\in P$ with $\{ \vec{k}\} = -\frac{k}r$. This form of
the blowup formula for Donaldson invariants and its higher rank analog
coincides with one given in \cite{MW,LNS1,MaM}.

\section{General gauge groups}\label{sec:gauge}

Our proof relies only on the blowup formula for degree $d=1,2$. Hence
it has a natural generalization to more general gauge groups. The
point is that we do not need the explicit formula \eqref{eq:Zsum} in
terms of Young tableaux.

Let $G$ be a compact semisimple Lie group. Let
$M^{\operatorname{reg}}(G,n)$ be the framed moduli space of
$G$-instantons on $S^4 = \R^4\cup\{ \infty\}$ with instanton number
$n$, which corresponds to $\pi_3(G) \cong \Z$. By \cite{AHS} it is a
nonsingular manifold, whose dimension can be computed by the index
theorem (and a standard calculation in the Lie algebra of $G$).
By the Hitchin-Kobayashi correspondence, the moduli space can be
identified with the framed moduli space of principal $G^c$-bundles on
$\CP^2 = \C^2\cup \ell_\infty$, where $G^c$ is the complexification of
$G$. When $G$ is a classical group, this version of the
Hitchin-Kobayashi correspondence was proved in \cite{Don} via
the ADHM description. Bando's analytic argument \cite{Bando} works for
arbitrary $G$. It is not clear, as far as the authors know,
whether we have a natural generalization of $M(r,n)$ for the group
$G$. Thus we can use only the Uhlenbeck compactification
\(
   M_0(G,n) =
   \bigsqcup_{m\le n} M^{\operatorname{reg}}(G,m)\times S^{n-m}\C^2.
\)
We also consider the framed moduli spaces
$\bM^{\operatorname{reg}}(G,k,n)$ and its Uhlenbeck compactification
$\bM_0(G,k,n)$ on the blowup. Here $k$ is the characteristic class in
$H^2$, which is considered as an element in $\pi_1(G)$.

Let $T$ be a maximal torus of $G$. Then we have an action of $\hT =
T^2\times T$ on the moduli spaces $M_0(G,n)$, $\bM_0(G,k,n)$.
Let $H^T_*(M_0(G,n))$, $H^T_*(\bM_0(G,k,n))$ denote the equivariant
homology groups.
The only fixed point in $M_0(G,n)$ is the ideal instanton consisting
of the trivial connection and the singularity concentrated at the
origin. We denote this point by $0$, and the inclusion $0\to M_0(G,n)$
by $\iota_0$. We {\it assume\/} that the localization theorem is
applicable to $M_0(G,n)$. This is guaranteed when $M_0(G,n)$ can be
equivariantly embedded in a finite dimensional representation of
$\hT$, or $M_0(G,n)$ can be endowed with a structure of
$\hT$-algebraic variety.  We define the {\it partition function\/} by
\begin{equation*}
   Z(\ve_1,\ve_2,\vec{a};\q)
   = \sum_{n=0}^\infty \q^n
   (\iota_{0*})^{-1} [M_0(G,n)],
\end{equation*}
where $[M_0(G,n)]$ is the fundamental class of $M_0(G,n)$. The
fundamental class is defined since the singular locus is lower
dimensional as fundamental classes of algebraic cycles are always defined.

\begin{Proposition}
The fixed points in $\bM_0(G,k,n)$ are parametrized by triples
$(\vec{k},l,m)$ where $\vec{k}\in\pi_1(T) \cong \Hom(S^1,T)$ and $l$,
$m$ are nonnegative integers. They satisfy the constraint
$\rho(\vec{k}) = k$ and $\frac12 (\vec{k},\vec{k}) + l + m = n$, where
$\rho$ is the homomorphism $\pi_1(T)\to \pi_1(G)$ induced by the
inclusion $T\subset G$, and $(\ ,\ )$ is the inner product on
$\operatorname{Lie}T$ such that the square of the length of the
highest root $\theta$ with respect to the induced inner product on the
dual space $\operatorname{Lie}T^*$ is equal to $2$.
\end{Proposition}

If we choose simple coroots $\alpha_i^\vee$ ($1\le i\le \dim T = \rank
G$), $\vec{k}$ can be identified with an $r$-tuple of rational numbers
$(k^1,\dots, k^r)\in\Q^r$ by $\vec{k} = \sum_i k^i\alpha_i^\vee$.

\begin{proof}
A fixed point in $\bM_0(G,k,n)$ is $(A,l[p_1] + m[p_2])$, where $A$ is
a reducible instanton (or a $G^c$-principal bundle which is reducible
to a $T^c$-bundle) with instanton number $n(A)$ and $l,m$ are integers
with $n(A) + l + m = n$. A reducible instanton $A$ on the blowup is
classified by $\vec{k}\in \pi_1(T) $. We have constraint
$\rho(\vec{k}) = k$, so that the induced bundle has the right
characteristic $k$. We also have
\begin{equation*}
   n(A) = \frac12 (\vec{k},\vec{k}),
\end{equation*}
where $(\ ,\ )$ is the inner product as above.
This can be proved as follows. Let $\mathfrak g^c$ be the
complexification of the vector bundle associated with the adjoint
representation. We have
\begin{equation*}
   c_2(\mathfrak g^c) = 
   \frac12 \sum_{\alpha\in\Delta} \langle\vec{k},\alpha\rangle^2
   = \frac12 B_G(\vec{k}, \vec{k})
   = h^\vee (\vec{k}, \vec{k}),
\end{equation*}
where $B_G$ is the Killing form, and $h^\vee$ is the dual Coxeter
number. For the last equality, see \cite[Exercise 6.1]{Kac}.
On the other hand, the instanton number is given by 
\(
   \frac{c_2(\mathfrak g^c)}{2h^\vee}.
\)
(See \cite[\S8]{AHS}.)
\end{proof}

For $G = \SU(r)$, the inner product $(\ ,\ )$ is the standard one used 
in earlier sections, and we have
\(
   h^\vee = r.
\)
Note that $c_2(\mathfrak g^c)$ is the complex dimension of the framed 
moduli space $\bM(G,k,n(A))$, so it is given by $2h^\vee n(A)$, as was 
shown in \cite{AHS}.

For a root $\alpha\in\Delta$, we define
\(
   l^{\vec{k}}_{\alpha}(\ve_1,\ve_2,\vec{a})
\)
by the same formula as as \eqref{eq:l}.
The Euler class of tangent space of $\bM^{\operatorname{reg}}(G,k,\frac12
(\vec{k},\vec{k}))$ at the reducible instanton $A$ is given by
\(
   \prod_{\alpha\in\Delta} l^{\vec{k}}_{\alpha}(\ve_1,\ve_2,\vec{a}).
\)

\begin{Conjecture}\label{conj1}
\textup{(1)} There exists a proper continuous map
$\widehat\pi_0\colon \bM_0(G,k,n)\to M_0(G,n')$ for some $n'$.

\textup{(2)} A neighborhood of the fixed point $(\vec{k},l,m)$ in
$\bM_0(G,k,n)$ is isomorphic to a neighborhood of
$(\vec{k},0,0)\times 0\times 0$ in $\bM_0(G,k,n-l-m)\times
M_0(G,l)\times M_0(G,m)$ as a $\hT$-space, where the $T^2$-actions on
the latter two factors are modified as $(t_1,t_2)\mapsto
(t_1,t_2/t_1)$ and $(t_1,t_2)\mapsto (t_1/t_2, t_2)$ respectively.
\end{Conjecture}

We define an equivariant cohomology class $\mu(C)\in
H^2_{\hT}(\bM^{\operatorname{reg}}(G,k,n))$ by
\begin{equation*}
   -\frac1{2 h^\vee} p_1(\widetilde{\mathfrak g}) / [C],
\end{equation*}
where $\widetilde{\mathfrak g}$ is the universal {\it adjoint\/}
bundle, i.e., the fiber is the Lie algebra $\mathfrak g$.

\begin{Conjecture}\label{conj2}
\textup{(1)} The class $\mu(C)$ extends to a class in
$H^2_{\hT}(\bM_0(G,k,n))$. We denote the extended class by the same
notation.

\textup{(2)} If $\iota_{(\vec{k},l,m)}$ denotes the inclusion of the
fixed point $(\vec{k},l,m)$ in $\bM_0(G,k,n)$, we have
\begin{equation*}
   \iota_{(\vec{k},l,m)}^{*}(\mu(C))
   = l\ve_1 + m\ve_2
   + (\vec{k},\vec{a})
   + \frac{(\vec{k},\vec{k})}2(\ve_1+\ve_2).
\end{equation*}
\end{Conjecture}

We define the {\it partition function\/} on the blowup by
\begin{equation*}
   \widehat{Z}(\ve_1,\ve_2,\vec{a};t;\q)
   = \sum_{n=0}^\infty \q^n \sum_{d=0}^\infty \frac{t^d}{d!}\;
   (\iota_{0*})^{-1}\widehat\pi_{0*}\left(\mu(C)^d\cap[{\bM_0(G,0,n)}]\right).
\end{equation*}
Then \eqref{eq:blowup} holds if we assume
Conjectures~\ref{conj1},~\ref{conj2}. \propref{prop:lowblowup} can be
modified as
\begin{equation*}
\begin{gathered}
   \widehat\pi_{0*}[\bM_0(G,0,n)] = [M_0(G,n)],
\\
   \widehat\pi_{0*}(\mu(C)^d\cap [\bM_0(G,0,n)]) = 0
   \quad\text{for $1\le d\le 2h^\vee-1$}.
\end{gathered}
\end{equation*}
The proof of the first equality is exactly the same. For the proof of
the second equality, we need a line bundle $\mathcal L$ and a section
which does not vanish on $\widehat\pi_0^{-1}(\overline{\{ 0\} \times
M_0^{\operatorname{reg}}(G,n-1)})$. I do not know such things exists
for genuine $\mu(C)$. But probably there exists such things for
$2h^\vee\mu(C)$. If this is indeed true, the rest of the argument is
the same as before.

We can now proceed as in the $\SU(r)$ case, we use this formula
$d=1,2$ to get (\ref{eq:ind1}, \ref{eq:ind2}). Considering the limit
$\ve_1,\ve_2\to 0$, we get \eqref{eq:Wh2} exactly as before. On the
other hand, the proof that the Seiberg-Witten prepotential satisfies
\eqref{eq:Wh2} was generalized to classical groups \cite{EGMM}.

\begin{Remark}
  The assumption that the localization theorem is applicable to
  $M_0(G,n)$ and $\bM(G,k,n)$ follows from the description in
  \cite{King} for a classical group $G$, since they are algebraic
  varieties. For general group, one can probably use the method in
  \cite{BFG}. Conjectures~\ref{conj1}, \ref{conj2} are true in view of
  King's description, except \ref{conj1}(2). We believe that
  \ref{conj1}(2) can be also checked, but we need a further study.
\end{Remark}

\end{document}